\documentclass[10pt]{article}
\usepackage{fancyhdr}
\usepackage{extramarks}
\usepackage{amsmath}
\usepackage{amsthm}
\usepackage{amsfonts}
\usepackage{siunitx}
\usepackage{tikz}
\usepackage[plain]{algorithm}
\usepackage{algpseudocode}
\usepackage{multirow}
\usepackage{booktabs}
\usepackage{palatino}
\usepackage{graphicx}
\usepackage{subfigure}
\usepackage[colorlinks,linkcolor=black,anchorcolor=black,citecolor=black,urlcolor=blue]{hyperref}
\usepackage{amsmath,bm}
\usepackage{booktabs}
\usepackage{mathtools}
\usepackage{amssymb}
\usepackage{caption}
\usepackage{capt-of}
\usepackage{mciteplus}
\usepackage{cite}
\usepackage{mathrsfs}
\usepackage{parskip}
\usepackage{soul}
\usepackage{textcomp}
\usepackage[colaction]{multicol}
\usepackage[switch]{lineno}
\usepackage{lipsum}
\usepackage{etoolbox}
\usepackage{longtable}
\usepackage{array}
\usepackage{tablefootnote}
\usepackage{ragged2e}
\usepackage{lscape}
\newcolumntype{C}[1]{>{\centering\arraybackslash}p{#1}}
\usetikzlibrary{automata,positioning}
\DeclareMathOperator{\im}{im}

\usetikzlibrary{automata,positioning}

\newtheorem{prop}{Proposition}[section]

\usepackage{tikz,xcolor,hyperref}
\definecolor{lime}{HTML}{A6CE39}
\DeclareRobustCommand{\orcidicon}{%
	\begin{tikzpicture}
	\draw[lime, fill=lime] (0,0) 
	circle [radius=0.16] 
	node[white] {{\fontfamily{qag}\selectfont \tiny ID}};
	\draw[white, fill=white] (-0.0625,0.095) 
	circle [radius=0.007];
	\end{tikzpicture}
	\hspace{-2mm}
}
\foreach \x in {A, ..., Z}{%
	\expandafter\xdef\csname orcid\x\endcsname{\noexpand\href{https://orcid.org/\csname orcidauthor\x\endcsname}{\noexpand\orcidicon}}
}

\begin{document}

\title{Persistent Path Laplacian}
\author{Rui Wang$^1$\footnote{Email: wangru25@msu.edu } and Guo-Wei Wei$^{1,2,3}$\footnote{
	Email: weig@msu.edu} \\
$^1$ Department of Mathematics, \\
Michigan State University, MI 48824, USA.\\
$^2$ Department of Electrical and Computer Engineering,\\
Michigan State University, MI 48824, USA. \\
$^3$ Department of Biochemistry and Molecular Biology,\\
Michigan State University, MI 48824, USA. \\
}
\date{\today} 

\maketitle

\begin{abstract}
 Path homology  proposed by S.-T.Yau and his co-workers provides a  new mathematical model for   directed graphs and networks.  Persistent path homology (PPH) extends the path homology with  filtration to deal with asymmetry structures. However, PPH is constrained to purely topological persistence and cannot track the homotopic shape evolution of data during filtration.  To overcome the limitation of PPH, persistent path Laplacian (PPL) is introduced  to capture the shape evolution of data. PPL's harmonic spectra fully recover PPH's topological persistence and its non-harmonic spectra  reveal the homotopic shape evolution of data during filtration. 
 
\end{abstract}
Keywords: Persistent path Laplacian, persistent path homology, topological Laplacian, spectral graph, directed graphs,  homotopic shape evolution of data.
%
 {\setcounter{tocdepth}{4} \tableofcontents}
 \newpage

\setcounter{page}{1}
\renewcommand{\thepage}{{\arabic{page}}}


\section{Introduction}
Recent years witness the emergence of a variety of advanced mathematical tools in topological data analysis (TDA) \cite{grbic2022aspects}. As the main workhorse of TDA, persistent homology (PH) \cite{carlsson2009topology,edelsbrunner2008persistent,zomorodian2005computing,edelsbrunner2000topological} pioneered a new branch in algebraic topology, offering a powerful tool to decode the topological structures of data during filtration in terms of persistent Betti numbers. Persistent homology has had tremendous  success in many areas of science and technology,  such as biology \cite{xia2014persistent}, chemistry \cite{townsend2020representation}, drug discovery \cite{nguyen2019mathematical}, 3D shape analysis \cite{skraba2010persistence}, etc. 

Inspired by the success of PH, other mathematical tools have been given due attention. One of them is de Rham-Hodge theory in differential geometry, which uses the differential forms to represent the cohomology of an oriented closed Riemannian manifold with boundary in terms of a topological Laplacian, namely Hodge Laplacian   \cite{dodziuk1977rham}.  The de Rham-Hodge theory has been applied to computational biology \cite{zhao2020rham}, graphic \cite{tong2003discrete}, and robotics \cite{mochizuki2009spatial}. However, like homology, the de Rham-Hodge theory does not offer an in-depth analysis of data, which is a famous problem in spectral geometry \cite{kac1966can}.  To overcome this drawback, the evolutionary de Rham-Hodge theory  \cite{chen2021evolutionary} was introduced in terms of persistent Hodge Laplacian to offer a multiscale analysis of the de Rham-Hodge theory. Defined on a family of evolutionary manifolds, the evolutionary de Rham-Hodge theory gives a new answer to, or at least reopens,  the famous 55-years old question ``can one hear the shape of a drum". \cite{kac1966can} The persistent Hodge Laplacian captures both the topological persistence and the homotopic shape evolution of data during filtration.

Nevertheless, the evolutionary de Rham-Hodge theory is set up on   Riemannian manifolds, which may be computationally demanding for large datasets. Hence, a similar multiscaled-based topological Laplacian,  called persistent spectral graph (PSG) \cite{wang2020persistent}, was proposed by introducing a filtration to combinatorial graph Laplacians. PSG, aka persistent Laplacian (PL) \cite{memoli2020persistent}, extends persistent homology to non-harmonic analysis of data, showing much advantage in sophisticated applications  \cite{meng2021persistent,chen2022persistent}. Dealing with point cloud data instead of manifolds, PL encodes a point cloud to a family of simplicial complexes generated from filtration and analyzes both harmonic and non-harmonic spectra. It is worthy to notice that the harmonic spectra from the null spaces of PLs reveal the same topological persistence like that of persistent homology, whereas, the non-harmonic spectra of PLs capture the homotopic shape evolution of data during the filtration. Meanwhile, open-source software called HERMES \cite{wang2021hermes} was developed for the simultaneous topological and geometric analysis of data.   However, like persistent homology, PSG treats all data points equally. That is to say, each point does not carry any labeled information such as the type, mass, color, etc. Therefore, an extension of PSG,   called persistent sheaf Laplacian  (PSL), was proposed to generalize cellular sheaves \cite{shepard1985cellular,hansen2019toward} for the multiscale analysis of point cloud data with attached labeled information \cite{wei2021persistent}. PSL is also a topological Laplacian that carries topological information in its null space but tracks homotopic shape evolution during filtration. Another interesting development is the persistent Dirac Laplacian (PDL) by  Ameneyro,  Maroulas, and Siopsis \cite{ameneyro2022quantum}. PDL offers an efficient quantum computation of persistent Betti numbers across different scales. These new approaches have great potentials to deal with complex data in science and engineering. 

It is noticed that the aforementioned homologies and topological Laplacians are insensitive to asymmetry or directed relations, which limits their representational power in encoding structures that have directional information. For example, in gene regulation data,  the directions of gene regulations are indicated by arrowheads or  perpendicular edges  in  systems biology \cite{long2008systems}. Therefore, a technique that can deal with  directed graphs (digraphs) is of vital importance to inferring gene regulation relationships. Notably, the path homology \cite{grigor2012homologies} proposed by Grigor’yan,  Lin,  Muranov, and  Yau  provides  a powerful tool to analyze datasets with asymmetric structures using the path complex.  Particular cases of homologies of digraphs and their path cohomology were also discussed   \cite{grigor2012homologies,grigor2015cohomology}.  The notion of path homology of digraphs has a richer mathematical structure than the earlier homology and Laplacian, opening new directions for both pure and applied mathematics.  For example,  path homology theory was extended  to various objects such as quivers, multigraphs, digraphs pairs, cylinder, cone, hypergraphs, etc. \cite{grigor2018path,grigor2018path2,grigor2019homology} Path homology has drawn much attention from researchers in the TDA community. To encode richer information, Chowdhury and   M\'{e}moli extended path homology to a persistent framework on a directed network \cite{chowdhury2018persistent}. Wang, Ren, and Wu constructed a weighted path homology for weight digraphs and proved a persistent version of a K\"{u}nneth-type formula for joins of weighted digraphs \cite{lin2019weighted}. Recently,   Dey,   Li, and   Wang have designed an efficient algorithm for $1$-dimensional persistent path homology \cite{dey2020efficient}, which is useful in real applications.  

Similar to persistent homology,  persistent path homology cannot track the homotopic shape evolution of data during filtration. To overcome this limitation, we introduce path Laplacian as a new topological Laplacian to analyze the spectral geometry of data, in addition to its topology.  
Moreover, we introduce a filtration to path Laplacian to obtain a persistent path Laplacian (PPL), a new framework that captures both the topological persistence and shape evolution of directed graphs and networks. By varying the filtration parameter, one can construct a series of digraphs, which result in a family of persistent path Laplacian matrices. The harmonic spectra of the persistent path Laplacian recover all the topological invariants of the digraphs, while the non-harmonic spectra provide additional geometric information, which can distinguish two systems when they are homotopy but geometrically different. PPL  has potential applications in science, engineering, industry, and technology. This work is organized as follows: Section 2 reviews the necessary background on path homology. Section 3 describes path Laplacian and persistent path Laplacian. Detailed PPL matrix constructions are illustrated with various examples for the interested readers { in Section 3 and Section 4.}

\section{Background on Path Homology and Directed Graph}\label{sec:Background}

Graph structure offers a powerful and versatile data representation that encodes inter-dependencies among constituents, which has been driven by widely spread applications in various fields such as graph theory, topological data analysis, science, and engineering. In this section, we first recap basic concepts in path homology, including paths on a finite set,  boundary operator on the path complex, and homologies of path complex. Then, we briefly review the concept of directed graphs (digraphs) and give a discussion of path homologies on the loopless directed graph. Such concepts and notations, due to Yau and coworkers, form a basis for us to introduce path Laplacian and persistent path Laplacian in \autoref{sec:ppl}.

\subsection{Paths on a Finite Set}
Denote set $V$  an arbitrary nonempty finite set, and elements in $V$ are called  vertices. For $p\in \mathbb{Z}^+_0$ { (i.e., a set with integers $p\ge 0$)}, an {\it elementary $p$-path} on $V$ is any sequence $i_0\dots i_p$ of $p+1$ vertices in $V$. An elementary $p$-path is an empty set $\emptyset$ for $p = -1$. For a fixed field $\mathbb{K}$, a   vector space that consists of all formal linear combinations of elementary $p$-paths with its coefficients in $\mathbb{K}$ is called {the space generated by the elementary paths}, denoted as $\Lambda_p = \Lambda_p(V, \mathbb{K}) = \Lambda_p(V)$.   One says the elements in $\Lambda_p$ are {\it $p$-paths} on $V$, and an elementary $p$-path $i_0\dots i_p  \in \Lambda_p$ is denoted by $e_{i_0\dots i_p}$. By definition, $\forall v\in \Lambda_p$, its unique representation can be given by the basis in $\Lambda_p$:
\begin{equation}
    v = \sum_{i_0, \dots, i_p \in V} c^{i_0\dots i_p}e_{i_0\dots i_p},
\end{equation}
where $c^{i_0\dots i_p}$ is the coefficient in $\mathbb{K}$. For instance, $\Lambda_0$ contains all linear combination of $e_i$ with $i\in V$, $\Lambda_1$ has all linear combination of $e_{ij}$ with $(i,j) \in V \times V$, and so on so forth. Since $\Lambda_{-1}$ consists of all multiples of $e$,   one has $\Lambda_{-1} \cong \mathbb{K}$. 

Additionally, $\forall p\in \mathbb{Z}^+_0$, the linear {\it boundary operator} from $\Lambda_p$ to $\Lambda_{p-1}$ that acts on elementary paths can be defined as 
\begin{equation}
    \partial: \Lambda_p \to \Lambda_{p-1}
\end{equation}
with
\begin{equation}\label{eq:boundary notation}
    \partial e_{i_0\dots i_p} = \sum_{q = 0}^{p} (-1)^q e_{i_0\dots \hat{i}_{q} \dots i_p},
\end{equation}
where $\hat{i}_{q}$ denotes the omission of index $i_{q}$ from the elementary $p$-path $e_{i_0\dots i_p}$.
One sets $\Lambda_{-2} = \{0\}$, and for $p = -1$, defines $\partial: \Lambda_{-1} \to \Lambda_{-2}$ to be a zero map. Following Lemma 2.1 in \cite{grigor2020path}, one has $\partial^2 = 0$, which indicates that the collection of boundary operator $\partial$ and space $\Lambda_p$ can form a chain complex of $V$ denoted as $\Lambda_{*} = \{\Lambda_p\}$ as 
\begin{equation}
    \cdots \Lambda_p \stackrel{\partial} \longrightarrow \Lambda_{p-1} \stackrel{\partial} \longrightarrow \cdots \stackrel{\partial} \longrightarrow  \Lambda_0 \stackrel{\partial} \longrightarrow \mathbb{K} \stackrel{\partial} \longrightarrow 0.
\end{equation}

Next, the concepts of regular path and non-regular path are introduced according to \cite{grigor2020path}. An elementary path $e_{i_0 \dots i_p}$ on a set $V$ is {\it regular} if $i_{k-1} \neq i_k$, and {\it non-regular} if $i_{k-1} = i_k$ for $k = 1,\dots, p$. For any $p\in \mathbb{Z}^{+}_{0} \cup \{-1\}$, let $\mathcal{R}_p$ be the subspace of $\Lambda_p$ spanned by all regular elementary paths, and $\mathcal{N}_p$ be the subspace of $\Lambda_p$ spanned by all non-regular elementary paths. Therefore, one has 
\begin{align*}
    \mathcal{R}_p &= \text{span}\{e_{i_0 \dots i_p}: i_0 \dots i_p \text{ is regular}\} \\
    \mathcal{N}_p &= \text{span}\{e_{i_0 \dots i_p}: i_0 \dots i_p \text{ is non-regular}\}.
\end{align*}
Note that $\mathcal{R}_p = \Lambda_p$ for integers { $p = -1, 0$}.

Then $\forall p\in \mathbb{Z}^{+}_{0} \cup \{-1\}, \ \Lambda_p = \mathcal{R}_p \oplus \mathcal{N}_p$. Therefore,
\begin{equation*}\label{eq:isomorphism}
    \mathcal{R}_p \cong \Lambda_p / \mathcal{N}_p.
\end{equation*}
According to Section 2.4 in \cite{grigor2020path}, the boundary operator $\partial$ is well-defined on the quotient space $\Lambda_p / \mathcal{N}_p$. Moreover, $\partial^2 = 0$ and the product rules are satisfied in the quotient space $\Lambda_p / \mathcal{N}_p$ as well.   One has an induced {\it regular boundary operator}:
\begin{equation}
    \bar{\partial}: \mathcal{R}_p \to \mathcal{R}_{p-1},
\end{equation}
where the regular boundary operator $ \bar{\partial}$ satisfies (\ref{eq:boundary notation}) except that all non-regular terms on the right hand side should be treated as $0$. Then a chain complex of $V$, denoted as $\mathcal{R}_{*}(V) = (\mathcal{R}_p)_p$ and equipped with $\bar{\partial}$,  can be expressed as:
\begin{equation}\label{eq:boundary for quotient space}
    \cdots \mathcal{R}_p \stackrel{\bar{\partial}} \longrightarrow \mathcal{R}_{p-1} \stackrel{\bar{\partial}} \longrightarrow \cdots \stackrel{\bar{\partial}} \longrightarrow  \mathcal{R}_0 \stackrel{\bar{\partial}} \longrightarrow \mathbb{K} \stackrel{\bar{\partial}} \longrightarrow 0.   
\end{equation}
It can be verified that $R_{p}\cong \Lambda_{p}/N_{p}$ is an isomorphism of chain complexes \cite{grigor2015cohomology}. In the following sections, for simplicity, we use $\partial$ to denote the boundary operator of Eq. (\ref{eq:boundary for quotient space}) unless specified differently.

\subsection{Path Complex}
A {\it path complex} over set $V$ is a nonempty collection $P$ of elementary paths on $V$ for any $n \in \mathbb{Z}^{+}_{0}$,
\begin{equation}\label{eq:path complex property}
    \text{if } i_0 \dots i_n \in P \text{, then } i_0 \dots i_{n-1} \in P, \text{ and } i_1 \dots i_n \in P.
\end{equation} 
For a fixed path complex, all the paths from $P$ are called {\it allowed} (i.e. $i_{k-1} \to i_k$ for any $k = 1,\dots,n$), while the elementary paths {on $V$} that are not in $P$ are {\it non-allowed}. We  say  a path complex $P$ is {\it perfect} if any subsequence of any path from $P$ is also in $P$. We choose $P_n$ to denote all $n$-paths from $P$. Then the set $P_{-1}$ has a single empty path $e$, the set $P_0$ consists of all the {\it vertices} of $P$, and clearly, $V = P_0$. To be noted, a path complex $P$ is a collection $\{P_n\}^{\infty}_{n=-1}$ satisfying Eq. (\ref{eq:path complex property}).
Let  $\mathcal{K}$ be an abstract simplicial complex defined over a finite vertex set $V$, satisfying
\begin{equation*}
    \text{if } \sigma \in \mathcal{K,} \text{ then any subset of } \sigma \text{ is also in } \mathcal{K}.
\end{equation*}
The collection of elementary paths on $V$ is denoted by $P(\mathcal{K})$. Follows from \cite{grigor2020path} (cf. Example 3.2), the family $P(\mathcal{K})$ is a path complex, and the allowed $n$-paths are $n$-simplices.

\subsection{Path Homology}\label{subsec:path homology}

For any $n\in \mathbb{Z}^{+}_0  $, the $\mathbb{K}$-linear space $\mathcal{A}_n$ is spanned by all the elementary $n$-paths from a given path complex $P = \{P_n\}^{\infty}_{n=0}$ over a finite set $V$, i.e.,
\begin{equation*}
    \mathcal{A}_n = \mathcal{A}_n(P) = \text{span} \{ e_{i_0 \dots i_n}: i_0 \dots i_n \in P_n\}.
\end{equation*} 
We call the elements of  $\mathcal{A}_n$ the $\it allowed \text{ } n\text{-}paths$. By the definition of {  $\mathcal{A}_n$, $\mathcal{A}_n\subset \Lambda_n$, and $\mathcal{A}_n = \Lambda_n$ for $n \leq 0$.} It is natural that the boundary operator $\partial$ defined on {  $\mathcal{R}_n$} can be introduced to $\mathcal{A}_n$ under certain condition: $\partial \mathcal{A}_n \subseteq \mathcal{A}_{n-1}$. For example, for perfect path complexes, we can obtain a chain complex:
\begin{equation*}
    \cdots \mathcal{A}_n \stackrel{\partial} \longrightarrow \mathcal{A}_{n-1} \stackrel{\partial} \longrightarrow \cdots \stackrel{\partial} \longrightarrow  \mathcal{A}_0 \stackrel{\partial} \longrightarrow \mathbb{K} \stackrel{\partial} \longrightarrow 0.
\end{equation*}

Next, we consider a general path complex $P$ (i.e., $\partial \mathcal{A}_n$ does not have to be a subspace of $\mathcal{A}_{n-1}$). For any $n\in \mathbb{Z}^{+}_0 \cup \{-1\} $, we define a subspace of $\mathcal{A}_{n}$:
\begin{equation}
    \Omega_n = \Omega_n(P) = \{v\in \mathcal{A}_n: \partial v \in \mathcal{A}_{n-1} \}.
\end{equation}
The elements of $\Omega_n$ are called $\partial$-invariant $n$-paths. To be noted, $\partial \Omega_n \subset \Omega_{n-1}$ always satisfies. Moreover, $\partial^2 = 0$ has been established in the previous section. Therefore, the $\it augmented$ chain complex of $\partial$-invariant paths can be denoted as
\begin{equation}\label{eq:augmented chain complex}
    \cdots \Omega_n \stackrel{\partial} \longrightarrow \Omega_{n-1} \stackrel{\partial} \longrightarrow \cdots \stackrel{\partial} \longrightarrow  \Omega_0 \stackrel{\partial} \longrightarrow \mathbb{K} \stackrel{\partial} \longrightarrow 0,
\end{equation}
whose homology group $\tilde{H}_n(P)$ of the chain complex in Eq. (\ref{eq:augmented chain complex}) are called the {\it reduced path homology groups} of the path complex $P$ for $n\in \mathbb{Z}^{+}_0 \cup \{-1\}$.
The truncated version of the chain complex in Eq. (\ref{eq:augmented chain complex}) for $n\in \mathbb{Z}^{+}_0$ is:
\begin{equation}\label{eq:standard chain complex}
    \cdots \Omega_n \stackrel{\partial} \longrightarrow \Omega_{n-1} \stackrel{\partial} \longrightarrow \cdots \stackrel{\partial} \longrightarrow  \Omega_0 \stackrel{\partial} \longrightarrow  0,
\end{equation}
whose homology group $H_n(P)$ of the chain complex in Eq. (\ref{eq:standard chain complex}) are called the {\it path homology groups} of the path complex $P$.

\subsection{Path Homology of Directed Graphs}\label{subsec:Path Homologies on Directed Graph}
A directed graph is an ordered pair $G = (V, E)$, where $V$ is a set of all vertices and $E$
is a set of ordered pairs of vertices (i.e. directed edges that satisfy $E \subseteq V \times V$). If $G = (V, E)$ does not contain any loop and multiple edge, then it is called  {\it simple directed graph}. Moreover, for the path homology of {\it multigraph} or {\it quiver}, one can refer to Ref. \cite{chartrand1977introductory}. In the following section of this work, we use $G(V,E)$ to represent the simple directed graphs unless specified differently. 

The path complex $P(G)$ is regular if $G=(V,E)$ is a simple directed graph. In this section, we mainly discuss the regular spaces $\Omega_n(G) = \Omega_n(P(G))$ and their associated regular homology groups $H(G) = H_n(P(G))$. Similar to the discussion in \autoref{subsec:path homology}, given a simple digraph $G(V,E)$, for any $n\in \mathbb{Z}^{+}_{0} \cup \{-1\}$, the space of $\partial$-invariant $n$-paths on $G$ is defined by the subspace of $\mathcal{A}_n(G) = \mathcal{A}_n(V, E; \mathbb{K})$:
\begin{equation*}
    \Omega_n = \Omega_n(G) = \{v\in \mathcal{A}_n: \partial v \in \mathcal{A}_{n-1} \},
\end{equation*}
with $\Omega_{-1} = \mathcal{A}_{-1} \cong \mathbb{K}$ and $\Omega_{-2} = \mathcal{A}_{-2} = \{0\}$. Since $\partial(\Omega_n) \subseteq \Omega_{n-1}$ (as $\partial^2 = 0$), then we have the following chain complex of $V$ denoted as $\Omega_{*}(V) = \{\Omega_n\}$,
\begin{equation*}\label{eq:digraph chain complex}
    \cdots \stackrel{\partial} \longrightarrow \Omega_3 \stackrel{\partial} \longrightarrow \Omega_{2} \stackrel{\partial} \longrightarrow \Omega_1 \stackrel{\partial} \longrightarrow \Omega_0 \stackrel{\partial} \longrightarrow \mathbb{K} \stackrel{\partial} \longrightarrow 0,
\end{equation*}
and the associated $n$-{\it dimensional path homology groups} of $G = (V,E)$ are defined as:
\begin{equation}
    H_n(G) = H_n(V,E;\mathbb{K}):= \ker(\partial |_{\Omega_n}) / \im (\partial|_{\Omega_{n+1}}).
\end{equation}
To be noted, the elements of $\ker(\partial |_{\Omega_n})$ are called {\it n-cycles}, and the elements of $\im (\partial|_{\Omega_{n+1}})$ are referred to as {\it n-boundaries}. For simplicity, we define $\partial_{n} = \partial|_{\Omega_n}$, and the chain complex of $\partial$-invariant paths is written as 
\begin{equation*}
    \cdots \Omega_{n+1} \stackrel{\partial_{n+1}} \longrightarrow \Omega_{n} \stackrel{\partial_n} \longrightarrow \Omega_{n-1} \stackrel{\partial_{n-1}} \longrightarrow  \Omega_{n-2} \cdots.
\end{equation*}
{ Notably, the path cohomology, introduced in Refs.  \cite{grigor2015cohomology,gomes2019path}, is isomorphic to the dual space of path homology when the coefficient ring is a field.}
The associated $n$-{\it dimensional path homology groups} of digraphs are defined as:
\begin{equation}
    H^n(G) = H^n(V,E;\mathbb{K}):= \ker(d_{n+1}) / \im (d_{n}),
\end{equation}
where $d$ is called coboundary operator. 

Given two simple digraphs $G = (V, E)$ and $G^{\prime} = (V^{\prime}, E^{\prime})$. According to the Definition 2.2 in \cite{grigor2014homotopy}, a {\it morphism of digraphs/digraphs map} from $G$ to $G^{\prime}$ is a map $f: V \to V^{\prime}$ such that for any directed edge $i \to j$ in $E$, one has either $f(i) \to f(j)$ is a directed edge on $E^{\prime}$ or $f(i) = f(j)$.

Let $f$ be a digraph map from $G$ to $G^{\prime}$. For $n\in \mathbb{Z}^{+}_{0}\cup \{-1\}$, one defines a map $(f_{**})_n: \Lambda_n(V) \to \Lambda_n(V^{\prime})$ such that:
\begin{equation}
    (f_{**})_n(e_{i_0 \dots i_n}) = e_{f(i_0) \dots f(i_n)}.
\end{equation}
Assume $\partial$ and $\partial^{\prime}$ are the boundary operators of chain complexes  $\Lambda_{*}(V)$ and $\Lambda_{*}(V^{\prime})$, then for $e_{i_0\dots i_n} \in \Lambda_n$, one has
\begin{align}\label{eq:chain map proof}
    ((f_{**})_{n-1} \circ \partial)(e_{i_0\dots i_n}) 
    &= \sum_{q=0}^n(-1)^q (f_{**})_{n-1}(e_{i_0\dots \hat{i}_{q} \dots i_n})\\
    &= \sum_{q=0}^n(-1)^q(e_{f(i_0) \dots \hat{f}(i_{q}) \dots f(i_n)})\\
    &= (\partial^{\prime}\circ (f_{*})_{n}) (e_{i_0\dots i_n}).
\end{align} 
Hence $f_{**}$ is a chain map.
By the definition of digraph map, $(f_{**})_n$ maps non-regular elementary $n$-paths on $V$ to non-regular elementary $n$-paths on $V^{\prime}$. Therefore, one has $(f_{**})_n(\mathcal{N}_n(V)) \subseteq \mathcal{N}_n(V^{\prime})$, and then $(f_{**})_n$ descends to a quotient homomorphism of chain complexes: 
\begin{equation}\label{eq:quotient homomorphism}
    (\tilde{f}_{**})_n: \Lambda_n(V)/\mathcal{N}_n(V) \to \Lambda_n(V^\prime)/\mathcal{N}_n(V^\prime).
\end{equation}
It can be verified that $R_{p}\cong \Lambda_{p}/N_{p}$ is an isomorphism of chain complexes \cite{grigor2015cohomology}, then the map in (\ref{eq:quotient homomorphism}) induces a morphism of chain complexes:
\begin{equation}
    (f_{*})_n: \mathcal{R}_n(V) \to \mathcal{R}_n(V^\prime).
\end{equation}
Since  $(f_{**})_n$ maps non-regular paths to non-regular, then similarly to what Eq. (\ref{eq:chain map proof}) shows, $(f_{*})_n$ is also a chain map that follows:
\begin{equation}
  (f_{*})_n(e_{i_0\dots i_n}) :=
    \begin{cases}
      e_{f(i_0) \dots f(i_n)} & \text{if $e_{f(i_0) \dots f(i_n)}$ is regular,}\\
      0 & \text{otherwise.}
    \end{cases}       
\end{equation}
Following the Theorem 2.10 in \cite{grigor2014homotopy}, the induced map $(f_{*})_n$ induces a morphism of chain complexes:
\begin{equation}
    (f_{*})_n : \Omega_{n}(G;\mathbb{K}) \to \Omega_{n}(G';\mathbb{K})
\end{equation}
and consequently induces a homomorphism between the path homology groups:
\begin{equation}
    (f_{*})_n : H_{n}(G;\mathbb{K}) \to H_{n}(G';\mathbb{K}),  \quad n\geq 0.
\end{equation}

\subsection{Homologies of Directed Subgraphs}
\begin{figure}[ht!]
	\includegraphics[width=0.8\textwidth]{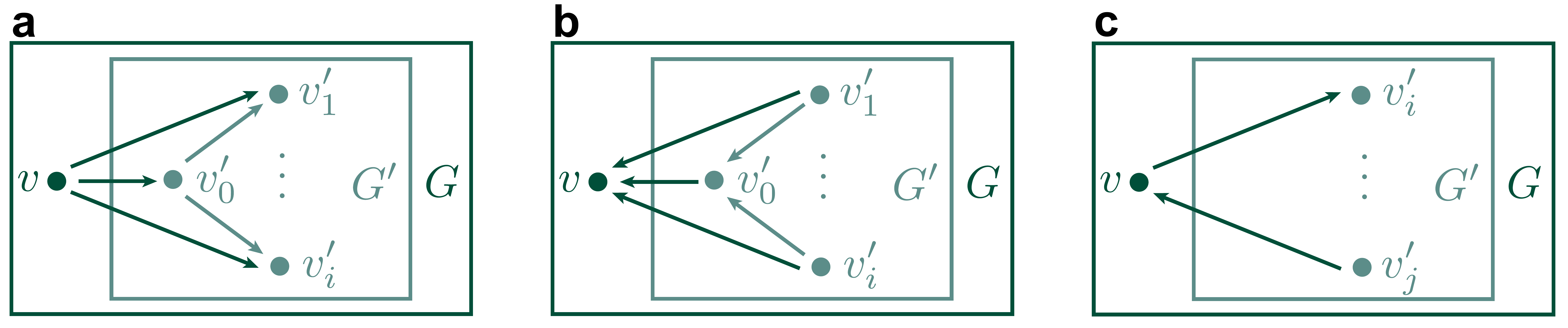}
	\centering
	\caption{Homologies of directed subgraphs. {\bf a}, {\bf b}, and {\bf c}  illustrate three subgraphs whose homology groups or homology group dimensions are related to the original digraphs.}
	\label{fig:homologies of subgraphs}
\end{figure}

 Some interesting propositions on the homologies of subgraphs provide a way to simplify complicated digraphs to relatively simple ones. Following the Section 4.2 in \cite{grigor2020path}, three propositions are discussed. 
\begin{prop}
    Given a simple digraph $G$ that has a vertex $v$ with $n$ outcoming arrows $v \to v^{\prime}_0, v \to v^{\prime}_1,\dots, v \to v^{\prime}_{n-1}$. Note that $v$ does not have any incoming arrows. Assume that for all $i \geq 1$, one has  $v^{\prime}_0 \to v^{\prime}_i$. Denote $G^{\prime}$ be the subgraph of $G$ by removing the vertex $v$ with all adjacent edges (i.e. $V^{\prime} = V\backslash \{v\}$ and $E^{\prime} = E\backslash \{v v^{\prime}_i\}_{i=0}^{n-1}$). Then, one has  $H_{*}(G) \cong H_{*}(G^{\prime})$ (See \autoref{fig:homologies of subgraphs} {\bf a}). 
\end{prop}

\begin{prop}
    Given a simple digraph $G = (V,E)$ that has a vertex $v$ with $n$ incoming arrows $v^{\prime}_0 \to v, v^{\prime}_1 \to v, \dots, v^{\prime}_{n-1} \to v$. Note that $v$ does not have any outcoming arrows. Assume that for all $i \geq 1$, one has $v^{\prime}_i \to v^{\prime}_0$. Denote $G^{\prime} = (V^{\prime}, E^{\prime})$ be the subgraph of $G$ by removing the vertex $v$ with all adjacent edges (i.e. $V^{\prime} = V\backslash \{v\}$ and $E^{\prime} = E\backslash \{v^{\prime}_i v\}_{i=0}^{n-1}$). Then, one has $H_{*}(G) \cong H_{*}(G^{\prime})$ (See \autoref{fig:homologies of subgraphs} {\bf b}). 
\end{prop}

\begin{prop}
    Given a simple digraph $G = (V,E)$ that has a vertex $v$ with only one outcoming arrow $v \to v^{\prime}_i$ and only one incoming arrow $v^{\prime}_j \to v$, where $i\neq j$. Denote $G^{\prime} = (V^{\prime}, E^{\prime})$ be the subgraph of $G$ (See \autoref{fig:homologies of subgraphs} {\bf c}) by removing the vertex $v$ and the adjacent edges $v \to v^{\prime}_i$ and $v^{\prime}_j \to v$ (i.e. $V^{\prime} = V\backslash \{v\}$ and $E^{\prime} = E\backslash \{v v^{\prime}_i,v^{\prime}_j v\}$). Then,
    \begin{enumerate}
        \item[(i)] $\dim H_p(G) = \dim H_p(G^{\prime})$ for $p \neq 2$ or for $p = 0,1$ if $v^{\prime}_j v^{\prime}_i$ is an edge/semi-edge in $G^{\prime}$.
        \item[(ii)] If $v^{\prime}_j v^{\prime}_i$ is neither an edge or a semi-edge in $G^{\prime}$, but $v^{\prime}_j$ and $v^{\prime}_i$ are in the same connected component of $G^{\prime}$, then $\dim H_1(G) = \dim H_1(G^{\prime} + 1)$, and $\dim H_0(G) = \dim H_0(G^{\prime})$.
        \item[(iii)] If $v^{\prime}_j$ and $v^{\prime}_i$ are not in the same connected component of $G^{\prime}$, then $\dim H_1(G) = \dim H_1(G^{\prime})$ and $\dim H_0(G) = \dim H_0(G^{\prime}) -1$.
    \end{enumerate}
\end{prop}


\section{Path Laplacian and Persistent Path Laplacian}\label{sec:ppl}

One can extract topological invariants by introducing the persistent Betti numbers from the homology groups along the filtration of simplicial complex \cite{zomorodian2005computing}. However, persistent Betti numbers do not capture homotopic geometric changes during filtration. Therefore, persistent topological Laplacians, including persistent Laplacian \cite{wang2020persistent,wang2021hermes} (persistent spectral graph) and persistent Hodge Laplacian \cite{chen2021evolutionary}, were introduced to reveal additional geometric information. Similarly, the constructions of path Laplacian and persistent path Laplacian are motivated by the earlier persistent spectral graphs \cite{wang2020persistent,wang2021hermes}. In this section, we first discuss the construction of path Laplacian. Then, we introduce   filtration to the path complex to generate a series of digraphs, which gives rise to persistent path Laplacian.

\subsection{Path Laplacian}\label{subsec:path laplacian}
Recall that a chain complex of $\partial$-invariant paths is given by
\begin{equation*}
    \cdots \Omega_{n+1} \stackrel{\partial_{n+1}} \longrightarrow \Omega_{n} \stackrel{\partial_n} \longrightarrow \Omega_{n-1} \stackrel{\partial_{n-1}} \longrightarrow  \Omega_{n-2} \cdots,
\end{equation*}
where $\Omega_n = \Omega_n(P) = \{v\in \mathcal{A}_n: \partial v \in \mathcal{A}_{n-1} \}$ and $\partial_n := \partial|_{\Omega_n}$. Alternatively,   assume $S_n:= S_n(P)$ to be the set of $n$-th elementary paths in $P$, then we define an inner product   
\begin{equation*}
    \langle \cdot, \cdot \rangle: S_n \times S_n \to \mathbb{R}
\end{equation*}
such that for any $e_{i_0\dots i_n}, e_{j_0\dots j_n} \in S_n$, the following satisfies
\begin{equation}
  \langle e_{i_0\dots i_n}, e_{j_0\dots j_n} \rangle =
    \begin{cases}
      1 & \text{if $e_{i_0\dots i_n} = e_{j_0\dots j_n}$,}\\
      0 & \text{otherwise.}
    \end{cases}       
\end{equation}

Let $M_{n}$ be a matrix representation of $\partial: \mathcal{A}_n \to \mathcal{A}_{n-1}$ with respect to the standard basis of $\mathcal{A}_{n}$ and $\mathcal{A}_{n-1}$. Define an inclusion map ${\displaystyle \iota_n :\Omega_n \hookrightarrow \mathcal{A}_n}$, then the matrix representation of ${\displaystyle \iota_n}$ with respect to the basis of $\Omega_n$ (i.e., the standard basis of $\mathcal{A}_n$ with the removal of {  generators} that are not in $\Omega_n$) and the standard basis of $\mathcal{A}_n$ is denoted as $O_n$. Denote the boundary matrix representation of $\partial_{n}$  as $B_{n}$, then we have
 
\begin{equation}\label{eq:boundary matrix 1}
    O_{n-1}B_{n} = \tilde{M}_{n} O_{n}. 
\end{equation}
If $O_{n-1}$ is a square matrix, then $O_n$ is actually an identity matrix, and we have 
\begin{equation}\label{eq:boundary matrix}
    B_{n} = O_{n-1}^{-1} \tilde{M}_{n} O_{n} = \tilde{M}_{n} O_{n},
\end{equation}
where $\tilde{M}_{n}$ is $M_{n}$ with the removal of rows that their basis are not elementary $(n-1)$-paths in $P$. Otherwise, $B_n$ is the least-square solution to Eq. (\ref{eq:boundary matrix 1}).
 
\begin{figure}[ht!]
	\includegraphics[width=0.3\textwidth]{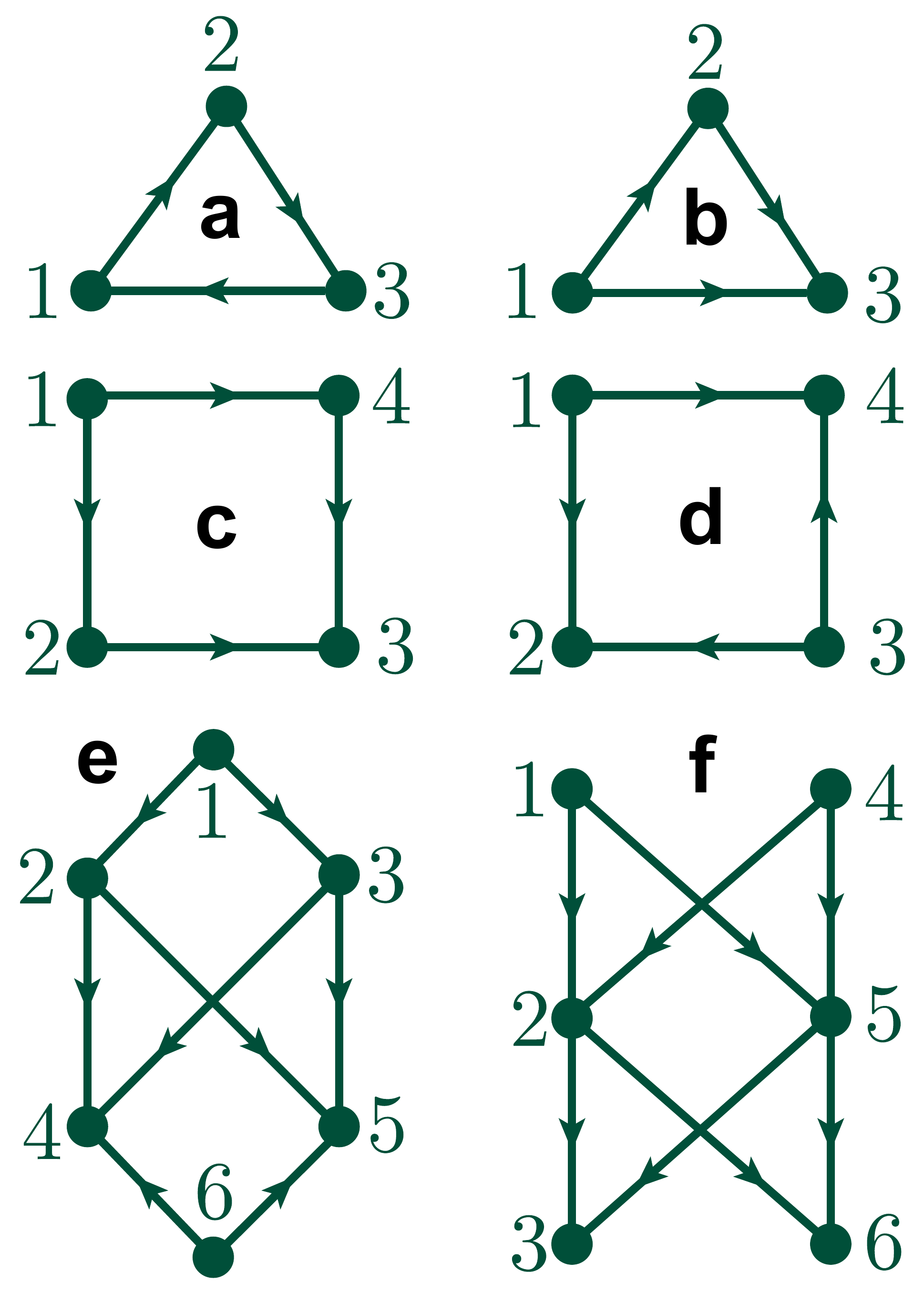}
	\centering
	\caption{Five digraphs. {\bf a} and {\bf b} Digraphs with 3 vertices and 3 directed edges. {\bf c} and {\bf d} Digraphs with 4 vertices and 4 directed edges. {\bf e} A digraph with 6 vertices and 8 directed edges. {\bf f} A digraph with 6 vertices and 8 directed edges.}
	\label{fig:path laplacians}
\end{figure}

 Note that $B_{n}$ is the matrix representation of $\partial_n$ with respect to the basis of $\Omega_n$ and $\Omega_{n-1}$. 
Dual space   $\Omega^{n}:= \text{Hom}(\Omega_n, \mathbb{K})$  of $\Omega_n$ is equipped with dual maps $d$ to form a cochain complex
\begin{equation*}
    \cdots \Omega^{n+1} \stackrel{d_{n+1}} \longleftarrow \Omega^{n} \stackrel{d_n} \longleftarrow \Omega^{n-1} \stackrel{d_{n-1}} \longleftarrow  \Omega^{n-2} \cdots,
\end{equation*}
where $d_{n}$ is called a coboundary operator. The inner product on $\Omega_n$ induces an inner product $\ll \cdot , \cdot \gg$ on $\Omega^n$ such that
\begin{equation*}
    \ll f , g \gg = \sum_{e\in S_n}f(e)g(e), \quad \forall f,g \in \Omega^n.
\end{equation*}
We denote the adjoint operator of $\partial_n$  be $\partial_n^*:  \Omega_{n-1} \to \Omega_n$. Note that similar inner product $\ll \cdot , \cdot \gg$ on $\Omega^n$ was  defined in  the literature \cite{horak2013spectra}. Hence, the coboundary operator $d_{n}$ is consistent with the adjoint operator $\partial_n^*$.
Then, { for integers $p\ge 0$}, the {\it $n$-th path Laplacian operator} is a linear operator: $\Delta_n: \Omega_n \to \Omega_n$ given by
\begin{equation}
    \Delta_n = \partial_{n+1} \partial_{n+1}^* + \partial_n^* \partial_n,
\end{equation}
and $\Delta_0 = \partial_{1} \partial_{1}^*$. The {\it $n$-th path Laplacian matrix} corresponding to $\Delta_n$ is expressed by
\begin{equation}\label{eq:path laplacian}
    L_n = B_{n+1} B_{n+1}^T + B_n^T B_n.
\end{equation}
Since  $L_n$ is positive semi-definite and symmetric, its eigenvalues are all real and non-negative. Additionally, recall that the Betti number $\beta_n$ of path complex $P$ satisfies
\begin{equation}
    \beta_n = \dim \ker \partial_n - \dim \im \partial_{n+1} = \dim \ker \Delta_n.
\end{equation}
It is easy to show that 
\begin{equation}
    \beta_n = \text{nullity} (L_n) = \text{the number of zero eigenvalues of } L_n.
\end{equation}
Moreover, assume the dimension of $L_n$ is $N$, then the set of spectra of $L_n$ is denoted as
\[
\text{Spectra}(L_n) = \{(\lambda_1)_n, (\lambda_2)_n, \cdots, (\lambda_N)_n  \}.
\]
\autoref{fig:path laplacians} shows 5 digraphs with multiple vertices and directed edges. Here, we take them as examples to give a detailed illustration of   $L_n$ matrix constructions.

\justifying
{\bf Construction of $L_0$ -- \autoref{fig:path laplacians}a} 
Since $L_0 = B_1B_1^T$, then we  first construct $B_1$, where $B_{1} = O_{0}^{-1} \tilde{M}_{1} O_{1}$ according to Eq. (\ref{eq:boundary matrix}), we have    
$O_0 = 
\bordermatrix{
        ~  &  e_{1} & e_{2} & e_{3}  \cr
        e_1 &  1  &  0  &  0 \cr
        e_2 &  0  &  1  &  0 \cr
        e_3 &  0  &  0  &  1 \cr
}$,
$M_1 = 
\bordermatrix{
        ~  &  e_{12} & e_{23} & e_{31}  \cr
        e_1 & -1  &  0  &  1 \cr
        e_2 &  1  & -1  &  0 \cr
        e_3 &  0  &  1  & -1 \cr
}$, and 
$O_1 = 
\bordermatrix{
        ~  &  e_{12} & e_{23} & e_{31}  \cr
        e_{12} &  1  &  0  &  0 \cr
        e_{23} &  0  &  1  &  0 \cr
        e_{31} &  0  &  0  &  1 \cr
}$. Since $e_1, e_2,$ and $ e_3$ are all elementary $0$-paths (vertices),  $M_1 = \tilde{M}_1$.  We have  
$B_1 = O_{0}^{-1} \tilde{M}_{1} O_{1} = 
\bordermatrix{
        ~  &  e_{12} & e_{23} & e_{31}  \cr
        e_1 & -1  &  0  &  1 \cr
        e_2 &  1  & -1  &  0 \cr
        e_3 &  0  &  1  & -1 \cr
}$. Then
$L_0 = B_1B_1^T = 
\begin{pmatrix}
        2  & -1  & -1 \\
       -1  &  2  & -1 \\
       -1  & -1  &  2 \\
\end{pmatrix}
$, which gives   $\text{Spectra}(L_0) = \{0,3,3\}$ and thus, one finally has $\beta_0 = 1$.

\justifying
{\bf Construction of $L_1$ -- \autoref{fig:path laplacians}a} 
We have $L_1 = B_2B_2^T + B_1^T B_1$, where $B_1$ has been formed, so we focus on the construction of $B_{2} = O_{1}^{-1} \tilde{M}_{2} O_{2}$ according to Eq. (\ref{eq:boundary matrix}). Since    
$O_1 = 
\bordermatrix{
        ~  &  e_{12} & e_{23} & e_{31}  \cr
        e_{12} &  1  &  0  &  0 \cr
        e_{23} &  0  &  1  &  0 \cr
        e_{31} &  0  &  0  &  1 \cr
}$,
$M_2 = 
\bordermatrix{
        ~  &  e_{123} & e_{231} & e_{312}  \cr
        e_{11} &  0  &  0  &  0 \cr
        e_{12} &  1  &  0  &  1 \cr
        e_{13} & -1  &  0  &  0 \cr
        e_{21} &  0  & -1  &  0 \cr
        e_{22} &  0  &  0  &  0 \cr
        e_{23} &  1  &  1  &  0 \cr
        e_{31} &  0  &  1  &  1 \cr
        e_{32} &  0  &  0  & -1 \cr
        e_{33} &  0  &  0  &  0 \cr
}$, and $O_2$ is a $3\times 0$ empty matrix since $\Omega_2 = \{0\}$. Therefore, 
$B_2 = O_{1}^{-1} \tilde{M}_{2} O_{2}$ is a $3\times 0$ empty matrix. Additionally, 
$ L_1 = B_2B_2^T + B_1^T B_1 =
\begin{pmatrix}
        2  & -1  & -1 \\
       -1  &  2  & -1 \\
       -1  & -1  &  2 \\
\end{pmatrix}
$, where $\text{Spectra}(L_1) = \{0,3,3\}$ and thus, one finally has  $\beta_1 = 1$.

\justifying
{\bf Construction of $L_2$ -- \autoref{fig:path laplacians}a} 
We have $L_2 = B_3B_3^T + B_2^T B_2$, where $B_2$ is an empty matrix. Hence, we focus on the construction of $B_{3} = O_{2}^{-1} \tilde{M}_{3} O_{3}$ according to Eq. (\ref{eq:boundary matrix}). We have $\mathcal{A}_2 = \text{span} \{e_{123}, e_{231}, e_{312}\}$ and $\mathcal{A}_1 = \text{span} \{e_{12}, e_{23}, e_{31}\}$. Note that $\partial_2(e_{123}) = e_{23} - e_{13} + e_{12}$ where $e_{13}$ is not in $\mathcal{A}_1$. Hence, $e_{123}$ is not in $\Omega_2$. The same conclusion can be deduced for $e_{231}$ and $e_{312}$. Therefore, we have $\Omega_2 = \{0\}$, and it is straightforward to get that $L_2$ is an empty matrix.   

\justifying
{\bf Construction of $L_0$ -- \autoref{fig:path laplacians}b} 
Since $L_0 = B_1B_1^T$, then we should first construct $B_1$, where $B_{1} = O_{0}^{-1} \tilde{M}_{1} O_{1}$ according to Eq. (\ref{eq:boundary matrix}). Since    
$O_0 = 
\bordermatrix{
        ~  &  e_{1} & e_{2} & e_{3}  \cr
        e_1 &  1  &  0  &  0 \cr
        e_2 &  0  &  1  &  0 \cr
        e_3 &  0  &  0  &  1 \cr
}$,
$M_1 = 
\bordermatrix{
        ~  &  e_{12} & e_{13} & e_{23}  \cr
        e_1 & -1  & -1  &  0 \cr
        e_2 &  1  &  0  & -1 \cr
        e_3 &  0  &  1  &  1 \cr
}$, and 
$O_1 = 
\bordermatrix{
        ~  &  e_{12} & e_{13} & e_{23}  \cr
        e_{12} &  1  &  0  &  0 \cr
        e_{13} &  0  &  1  &  0 \cr
        e_{23} &  0  &  0  &  1 \cr
}$. Since $e_1, e_2, $ and $ e_3$ are all elementary $0$-paths (vertices). Therefore, $M_1 = \tilde{M}_1$, and we have  
$B_1 = O_{0}^{-1} \tilde{M}_{1} O_{1} = 
\bordermatrix{
        ~  &  e_{12} & e_{13} & e_{23}  \cr
        e_1 & -1  & -1  &  0 \cr
        e_2 &  1  &  0  & -1 \cr
        e_3 &  0  &  1  &  1 \cr
}$. Then
$L_0 = B_1B_1^T = 
\begin{pmatrix}
        2  & -1  & -1 \\
       -1  &  2  & -1 \\
       -1  & -1  &  2 \\
\end{pmatrix}
$, which gives the $\text{Spectra}(L_0) = \{0,3,3\}$ and thus, one finally has  $\beta_0 = 1$.

\justifying
{\bf Construction of $L_1$ -- \autoref{fig:path laplacians}b} 
We have $L_1 = B_2B_2^T + B_1^T B_1$, where $B_1$ has been formed, so we focus on the construction of $B_{2} = O_{1}^{-1} \tilde{M}_{2} O_{2}$ according to Eq.  (\ref{eq:boundary matrix}). First, $\mathcal{A}_2 = \text{span} \{e_{123}\}$ and $\mathcal{A}_1 = \text{span} \{e_{12}, e_{13}, e_{23}\}$. Note that $\partial_2(e_{123}) = e_{23} - e_{13} + e_{12}$ where $e_{12}, e_{23}, $ and $ e_{13}$ are all in $\mathcal{A}_1$. Hence, $\Omega_2 = \mathcal{A}_2 = \text{span}\{e_{123}\}$. Note that
$O_1 = 
\bordermatrix{
        ~  &  e_{12} & e_{13} & e_{23}  \cr
        e_{12} &  1  &  0  &  0 \cr
        e_{13} &  0  &  1  &  0 \cr
        e_{23} &  0  &  0  &  1 \cr
}$,
$M_2 = 
\bordermatrix{
        ~  &  e_{123}   \cr
        e_{11} &  0  \cr
        e_{12} &  1  \cr
        e_{13} & -1  \cr
        e_{21} &  0  \cr
        e_{22} &  0  \cr
        e_{23} &  1  \cr
        e_{31} &  0  \cr
        e_{32} &  0  \cr
        e_{33} &  0  \cr
}$, and 
$O_2 = 
\bordermatrix{
        ~  &  e_{123} \cr
        e_{123} &  1  \cr
}$. The $e_{11},e_{21},e_{22},e_{31},e_{32}, $ and $ e_{33}$ are not elementary $1$-paths in $P$. Hence, 
$\tilde{M}_2 = 
\bordermatrix{
        ~  &  e_{123}   \cr
        e_{12} &  1  \cr
        e_{13} & -1  \cr
        e_{23} &  1  \cr
}$, and then 
$B_2 = O_{1}^{-1} \tilde{M}_{2} O_{2} =
\bordermatrix{
        ~  &  e_{123}   \cr
        e_{12} &  1  \cr
        e_{13} & -1  \cr
        e_{23} &  1  \cr
}$.
Therefore,
$ L_1 = B_2B_2^T + B_1^T B_1 =
\begin{pmatrix}
        3  &  0  &  0 \\
        0  &  3  &  0 \\
        0  &  0  &  3 \\
\end{pmatrix}
$, where $\text{Spectra}(L_1) = \{3,3,3\}$ and thus, we finally have  $\beta_1 = 0$.

\justifying
{\bf Construction of $L_2$ -- \autoref{fig:path laplacians}b} 
According to Eq. (\ref{eq:path laplacian}), we have $L_2 = B_3B_3^T + B_2^T B_2$ and $B_{3} = O_{2}^{-1} \tilde{M}_{3} O_{3}$. Since there is no $3$-path existing, so the $M_3$ and $O_3$ are both empty matrix. Hence $L_2 = (3)$, $\text{Spectra}(L_2) = \{3\}$, and thus, one has $\beta_2 = 0$.

In the following section, we will omit the detailed construction steps of boundary matrix $B_n$. \autoref{tab:square case 1}, \autoref{tab:square case 2}, \autoref{tab:Pentagon}, and \autoref{tab:network} list the boundary matrix $B_n$ and the $n$-th path Laplacian matrix $L_n$ for with its corresponding Betti numbers $\beta_n$ and spectrum $\text{Spectra}(L_n)$ for \autoref{fig:path laplacians} {\bf c}, {\bf d}, {\bf e}, and {\bf f}. It is worth to mention that $\beta_n$ can distinguish the same graph with different paths assigned. For example, \autoref{fig:path laplacians} {\bf c} and {\bf d} have the same undirected graph structure with different paths assigned. We have $\beta_1=0$ for \autoref{fig:path laplacians} {\bf c} and  $\beta_1=1$ for  \autoref{fig:path laplacians} {\bf d}.

\begin{table}[H]
    \centering
    \setlength\tabcolsep{6pt}
    \captionsetup{margin=0.9cm}
    \caption{Illustration of digraph {\bf c} in \autoref{fig:path laplacians}}
    \begin{tabular}{c|ccc}
    \hline
    $n$  & $n=0$ & $n=1$ & $n=2$ \\ \hline 
    $\Omega_n$ & $\text{span}\{e_1,e_2,e_3,e_4\}$ & $\text{span}\{e_{12}, e_{14}, e_{23}, e_{43}\}$ & $\text{span}\{e_{143} - e_{123}\}$ \\ \hline \\
    $B_{n+1}$  & $\begin{array}{@{}r@{}c@{}c@{}c@{}c@{}c@{}l@{}}
            & e_{12} & e_{14} & e_{23} & e_{43}  \\
           \left.\begin{array}{c}
            e_1 \\
            e_2 \\
            e_3 \\
            e_4
            \end{array}\right(
            & \begin{array}{c} -1 \\  1  \\  0  \\  0  \end{array}
            & \begin{array}{c} -1 \\  0  \\  0  \\  1  \end{array}
            & \begin{array}{c}  0 \\ -1  \\  1  \\  0  \end{array}
            & \begin{array}{c}  0 \\  0  \\  1  \\ -1  \end{array}
            & \left)\begin{array}{c} \\ \\  \\ \\ 
            \end{array}\right.
        \end{array}$
        & $\begin{array}{@{}r@{}c@{}c@{}c@{}c@{}c@{}l@{}}
            & e_{143} - e_{123}   \\
           \left.\begin{array}{c}
            e_{12} \\
            e_{14} \\
            e_{23} \\
            e_{43} 
            \end{array}\right(
            & \begin{array}{c} -1 \\  1  \\  -1  \\  1  \end{array}
            & \left)\begin{array}{c} \\ \\  \\ \\ 
            \end{array}\right.
        \end{array}$ &  $1\times 0$ empty matrix    \\  \\
    $L_n$  & $\left(\begin{array}{ccccc}
         2  &  -1  &   0  &  -1  \\
        -1  &   2  &  -1  &   0  \\
         0  &  -1  &   2  &  -1  \\
        -1  &   0  &  -1  &   2   
    \end{array}\right)$  & $\left(\begin{array}{ccccc}
         3  &   0  &   0  &  -1  \\
         0  &   3  &  -1  &   0  \\
         0  &  -1  &   3  &   0  \\
        -1  &   0  &   0  &   3   
    \end{array}\right)$ &    $\left(\begin{array}{ccccc}
         4 
    \end{array}\right)$  \\  \\
    $\beta_{n}$                        & 1               & 0 & 0    \\  \\
    $\text{Spectra}(L_{n})$    & $\{0,2,2,4\}$   & $\{2,2,4,4\}$ & \{4\}  \\ \hline
    \end{tabular}
    \\
    \label{tab:square case 1}
\end{table}

\begin{table}[H]
    \centering
    \setlength\tabcolsep{6pt}
    \captionsetup{margin=0.9cm}
    \caption{Illustration of digraph {\bf d} in \autoref{fig:path laplacians}}
    \begin{tabular}{c|ccc}
    \hline
    $n$  & $n=0$ & $n=1$ & $n=2$ \\ \hline 
    $\Omega_n$ & $\text{span}\{e_1,e_2,e_3,e_4\}$ & $\text{span}\{e_{12}, e_{14}, e_{32}, e_{34}\}$ & $\{0\}$ \\ \hline \\
    $B_{n+1}$  & $\begin{array}{@{}r@{}c@{}c@{}c@{}c@{}c@{}l@{}}
            & e_{12} & e_{14} & e_{32} & e_{34}  \\
           \left.\begin{array}{c}
            e_1 \\
            e_2 \\
            e_3 \\
            e_4
            \end{array}\right(
            & \begin{array}{c} -1 \\  1  \\  0  \\  0  \end{array}
            & \begin{array}{c} -1 \\  0  \\  0  \\  1  \end{array}
            & \begin{array}{c}  0 \\  1  \\ -1  \\  0  \end{array}
            & \begin{array}{c}  0 \\  0  \\ -1  \\  1  \end{array}
            & \left)\begin{array}{c} \\ \\  \\ \\ 
            \end{array}\right.
        \end{array}$
        & $4\times 0$ empty matrix &  $\left(\begin{array}{ccccc}
         /
    \end{array}\right)$   \\  \\
    $L_n$  & $\left(\begin{array}{ccccc}
         2  &  -1  &   0  &  -1  \\
        -1  &   2  &  -1  &   0  \\
         0  &  -1  &   2  &  -1  \\
        -1  &   0  &  -1  &   2   
    \end{array}\right)$  & $\left(\begin{array}{ccccc}
         2  &   1  &   1  &   0  \\
         1  &   2  &   0  &   1  \\
         1  &   0  &   2  &   1  \\
         0  &   1  &   1  &   2   
    \end{array}\right)$ &    $\left(\begin{array}{ccccc}
         / 
    \end{array}\right)$  \\  \\
    $\beta_{n}$                        & 1               & 1 & 0    \\  \\
    $\text{Spectra}(L_{n})$    & $\{0,2,2,4\}$   & $\{0,2,4,4\}$ & / \\ \hline
    \end{tabular}
    \\
    \label{tab:square case 2}
\end{table}

\begin{table}[H]\footnotesize
    \centering
    \setlength\tabcolsep{0pt}
    \captionsetup{margin=0.9cm}
    \caption{Illustration of digraph {\bf e} in \autoref{fig:path laplacians}.}
    \begin{tabular}{c|ccc}
    \hline
    $n$  & $n=0$ & $n=1$ & $n=2$ \\ \hline 
    $\Omega_n$ & $\text{span}\{e_1,e_2,e_3,e_4,e_5,e_6\}$ & $\text{span}\{e_{12},e_{13}, e_{24}, e_{25},e_{34}, e_{35},e_{64},e_{65}\}$ & $\text{span}\{e_{134}-e_{124}, e_{135}-e_{125}\}$ \\\hline \\
    $B_{n+1}$  & $\begin{array}{@{}r@{}c@{}c@{}c@{}c@{}c@{}c@{}c@{}c@{}l@{}}
            & e_{12} & e_{13} & e_{24} & e_{25} & e_{34} & e_{35} & e_{64} & e_{65}  \\
           \left.\begin{array}{c}
            e_1 \\
            e_2 \\
            e_3 \\
            e_4 \\
            e_5 \\
            e_6
            \end{array}\right(
            & \begin{array}{c} -1 \\  1  \\  0  \\  0 \\ 0 \\ 0  \end{array}
            & \begin{array}{c} -1 \\  0  \\  1  \\  0 \\ 0 \\ 0  \end{array}
            & \begin{array}{c}  0 \\ -1  \\  0  \\  1 \\ 0 \\ 0  \end{array}
            & \begin{array}{c}  0 \\ -1  \\  0  \\  0 \\ 1 \\ 0  \end{array}
            & \begin{array}{c}  0 \\  0  \\ -1  \\  1 \\ 0 \\ 0  \end{array}
            & \begin{array}{c}  0 \\  0  \\ -1  \\  0 \\ 1 \\ 0  \end{array}
            & \begin{array}{c}  0 \\  0  \\  0  \\  1 \\ 0 \\-1  \end{array}
            & \begin{array}{c}  0 \\  0  \\  0  \\  0 \\ 1 \\-1  \end{array}
            & \left)\begin{array}{c} \\ \\  \\ \\ \\ \\
            \end{array}\right.
        \end{array}$
        & $\begin{array}{@{}r@{}c@{}c@{}c@{}c@{}c@{}l@{}}
            & e_{134}- e_{124}  & {\quad }e_{135} - e_{125}  \\
           \left.\begin{array}{c}
            e_{12} \\
            e_{13} \\
            e_{24} \\
            e_{25} \\
            e_{34} \\
            e_{35} \\
            e_{64} \\
            e_{65} 
            \end{array}\right(
            & \begin{array}{c} -1 \\  1  \\ -1  \\  0 \\ 1 \\ 0 \\ 0 \\ 0 \end{array}
            & \begin{array}{c} -1 \\  1  \\  0  \\ -1 \\ 0 \\ 1 \\ 0 \\ 0 \end{array}
            & \left)\begin{array}{c} \\ \\  \\ \\ \\ \\ \\ \\
            \end{array}\right.
        \end{array}$ &  $2\times 0$ empty matrix    \\  \\
    $L_n$  & $\left(\begin{array}{cccccccccccc}
         2  &  -1  &  -1   &   0   &  0   &  0   \\
        -1  &   3  &   0   &  -1   & -1   &  0   \\
        -1  &   0  &   3   &  -1   & -1   &  0   \\
         0  &  -1  &  -1   &   3   &  0   & -1   \\
         0  &  -1  &  -1   &   0   &  3   & -1   \\
         0  &   0  &   0   &  -1   & -1   &  2   
    \end{array}\right)$  & $\left(\begin{array}{cccccccccccc}
         4  &  -1  &   0   &   0   & -1   & -1  &  0   &  0   \\
        -1  &   4  &  -1   &  -1   &  0   &  0  &  0   &  0   \\
         0  &  -1  &   3   &   1   &  0   &  0  &  1   &  0   \\
         0  &  -1  &   1   &   3   &  0   &  0  &  0   &  1   \\
        -1  &   0  &   0   &   0   &  3   &  1  &  1   &  0   \\
        -1  &   0  &   0   &   0   &  1   &  3  &  0   &  1   \\    
         0  &   0  &   1   &   0   &  1   &  0  &  2   &  1   \\    
         0  &   0  &   0   &   1   &  0   &  1  &  1   &  2   
    \end{array}\right)$ &    $\left(\begin{array}{cccccccccccc}
         4  &  2    \\
         2  &  4
    \end{array}\right)$ \\  \\
    $\beta_{n}$                        & 1               & 1 & 0    \\  \\
    $\text{Spectra}(L_{n})$    & $\{0,1.4384,3,3,3,5\}$   & $\{0,1.4384,2,3,3,3,5.5616,6\}$ & \{2,6\} \\ \hline
    \end{tabular}
    \\
    \label{tab:Pentagon}
\end{table}

\begin{table}[H]\scriptsize
    \centering
    \setlength\tabcolsep{0pt}
    \captionsetup{margin=0.9cm}
    \caption{Illustration of digraph {\bf f} in \autoref{fig:path laplacians}.}
    \begin{tabular}{c|ccc}
    \hline
    $n$  & $n=0$ & $n=1$ & $n=2$ \\ \hline 
    $\Omega_n$ & $\text{span}\{e_1,e_2,e_3,e_4,e_5,e_6\}$ & $\text{span}\{e_{12},e_{15}, e_{23}, e_{26},e_{42}, e_{45},e_{53},e_{56}\}$ & $\text{span}\{e_{153}-e_{123},$ \\
     &  &  & $e_{156}-e_{126},$ \\ 
     &  &  & $e_{453}-e_{423},$ \\ 
     &  &  & $e_{456}-e_{426}\}$\\\hline \\
    $B_{n+1}$  & $\begin{array}{@{}r@{}c@{}c@{}c@{}c@{}c@{}c@{}c@{}c@{}l@{}}
            & e_{12} & e_{15} & e_{23} & e_{26} & e_{42} & e_{45} & e_{53} & e_{56}  \\
           \left.\begin{array}{c}
            e_1 \\
            e_2 \\
            e_3 \\
            e_4 \\
            e_5 \\
            e_6
            \end{array}\right(
            & \begin{array}{c} -1 \\  1  \\  0  \\  0 \\ 0 \\ 0  \end{array}
            & \begin{array}{c} -1 \\  0  \\  0  \\  0 \\ 1 \\ 0  \end{array}
            & \begin{array}{c}  0 \\ -1  \\  1  \\  0 \\ 0 \\ 0  \end{array}
            & \begin{array}{c}  0 \\ -1  \\  0  \\  0 \\ 0 \\ 1  \end{array}
            & \begin{array}{c}  0 \\  1  \\  0  \\ -1 \\ 0 \\ 0  \end{array}
            & \begin{array}{c}  0 \\  0  \\  0  \\ -1 \\ 1 \\ 0  \end{array}
            & \begin{array}{c}  0 \\  0  \\  1  \\  0 \\-1 \\ 0  \end{array}
            & \begin{array}{c}  0 \\  0  \\  0  \\  0 \\-1 \\ 1  \end{array}
            & \left)\begin{array}{c} \\ \\  \\ \\ \\ \\
            \end{array}\right.
        \end{array}$
        & $\begin{array}{@{}r@{}c@{}c@{}c@{}c@{}c@{}l@{}}
            & e_{153}-e_{123} & {\quad }e_{156}-e_{126} & {\quad }e_{453}-e_{423} & {\quad }e_{456}-e_{426}  \\
           \left.\begin{array}{c}
            e_{12} \\
            e_{15} \\
            e_{23} \\
            e_{26} \\
            e_{42} \\
            e_{45} \\
            e_{53} \\
            e_{56} 
            \end{array}\right(
            & \begin{array}{c} -1 \\  1  \\ -1  \\  0 \\ 0 \\ 0 \\ 1 \\ 0 \end{array}
            & \begin{array}{c} -1 \\  1  \\  0  \\ -1 \\ 0 \\ 0 \\ 0 \\ 1 \end{array}
            & \begin{array}{c}  0 \\  0  \\ -1  \\  0 \\-1 \\ 1 \\ 1 \\ 0 \end{array}
            & \begin{array}{c}  0 \\  0  \\  0  \\ -1 \\-1 \\ 1 \\ 0 \\ 1 \end{array}
            & \left)\begin{array}{c} \\ \\  \\ \\ \\ \\ \\ \\
            \end{array}\right.
        \end{array}$ &  $4\times 0$ empty matrix    \\  \\
    $L_n$  & $\left(\begin{array}{cccccccccccc}
         2  &  -1  &   0   &   0   & -1   &  0   \\
        -1  &   4  &  -1   &  -1   &  0   & -1   \\
         0  &  -1  &   2   &   0   & -1   &  0   \\
         0  &  -1  &   0   &   2   & -1   &  0   \\
        -1  &   0  &  -1   &  -1   &  4   & -1   \\
         0  &  -1  &   0   &   0   & -1   &  2   
    \end{array}\right)$  & $\left(\begin{array}{cccccccccccc}
         4  &  -1  &   0   &   0   &  1   &  0  & -1   & -1   \\
        -1  &   4  &  -1   &  -1   &  0   &  1  &  0   &  0   \\
         0  &  -1  &   4   &   1   &  0   & -1  & -1   &  0   \\
         0  &  -1  &   1   &   4   &  0   & -1  &  0   & -1   \\
         1  &   0  &   0   &   0   &  4   & -1  & -1   & -1   \\
         0  &   1  &  -1   &  -1   & -1   &  4  &  0   &  0   \\    
        -1  &   0  &  -1   &   0   & -1   &  0  &  4   &  1   \\    
        -1  &   0  &   0   &  -1   & -1   &  0  &  1   &  4   
    \end{array}\right)$ &    $\left(\begin{array}{cccccccccccc}
         4  &  2  &  2  & 0  \\
         2  &  4  &  0  & 2  \\
         2  &  0  &  4  & 2  \\
         0  &  2  &  2  & 4  
    \end{array}\right)$ \\  \\
    $\beta_{n}$                        & 1               & 0 & 1    \\  \\
    $\text{Spectra}(L_{n})$    & $\{0,2,2,2,4,6\}$   & $\{2,2,2,4,4,4,6,8\}$ & \{0,4,4,8\} \\ \hline
    \end{tabular}
    \\
    \label{tab:network}
\end{table}

\subsection{Persistent Path Laplacian}

From Section \ref{subsec:path laplacian}, the way to calculate both harmonic spectra (topological invariants) and non-harmonic spectra of $n$-th path Laplacian matrix is genuinely free of metrics or coordinates, which contains too little information to fully describe the object. Therefore, inspired by the idea of the persistent spectral graph (PSG), persistent path Laplacian (PPL) is proposed to create a sequence of digraphs induced by varying a filtration parameter to encode more geometric or structural information.



First, we consider a {\it filtration of digraphs}  $\mathcal{G}: \mathbb{R} \to \mathcal{D}$, which is a morphism   $f_{s,t}: H_p(\mathcal{G}_t;\mathbb{K}) \to H_p(\mathcal{G}_s;\mathbb{K})$  from the category of real number $\mathbb{R}$ to the category of digraphs $\mathcal{D}$ that satisfies:
\begin{equation*}
    \mathcal{G}(t) \subseteq \mathcal{G}(s), \forall  t \leq s,
\end{equation*}
where $G_t:=\mathcal{G}(t) \in \mathcal{D}$ and $G_s:=\mathcal{G}(s)  \in \mathcal{D}$. Consider a { sequence of finitely many positive integers} $1,2,\dots,m$, we have a sequence of digraphs $$G_1 \subseteq G_2 \subseteq \dots \subseteq G_m.$$
For each digraph $G_t$, we denote its corresponding chain group to be $\Omega_n(G_t)$, and the $n$-boundary operator of $G_t$ is denoted by $\partial_n^t: \Omega_n(G_t) \to \Omega_{n-1}(G_t),  \forall  n \ge 0$ . 


Similarly, as in  persistent homology, a sequence of chain complexes can be denoted as 
\begin{equation}
    \left.\begin{array}{cccccccccccccc}
        \cdots & \Omega_{n+1}^1 &
        \xrightarrow[]{\partial_{n+1}^1} & \Omega_n^1 &
        \xrightarrow[]{\partial_n^1} & \cdots & \xrightarrow[]{\partial_3^1} & \Omega_2^1 & \xrightarrow[]{\partial_2^1} & \Omega_1^1 & \xrightarrow[]{\partial_1^1} & \Omega_0^1 & \xrightarrow[]{\partial_0^1} & \Omega_{-1}^1 \\
        & \rotatebox{-90}{$\hookrightarrow$} &  & \rotatebox{-90}{$\hookrightarrow$} &  &  &  & \rotatebox{-90}{$\hookrightarrow$} &  & \rotatebox{-90}{$\hookrightarrow$} &  & \rotatebox{-90}{$\hookrightarrow$} &  &  \\
        \cdots & \Omega_{n+1}^2 &
        \xrightarrow[]{\partial_{q+1}^2} & \Omega_n^2 &
        \xrightarrow[]{\partial_n^2} & \cdots &
        \xrightarrow[]{\partial_3^2} & \Omega_2^2 & \xrightarrow[]{\partial_2^2} & \Omega_1^2 & \xrightarrow[]{\partial_1^2} & \Omega_0^2 & \xrightarrow[]{\partial_0^2} & \Omega_{-1}^2 \\
				     & \rotatebox{-90}{  $\cdots$} &  & \rotatebox{-90}{  $\cdots$} &  &  &  & \rotatebox{-90}{ $\cdots$} &  & \rotatebox{-90}{  $\cdots$} &  & \rotatebox{-90}{ $\cdots$} &  &  \\
        & \rotatebox{-90}{$\hookrightarrow$} &  & \rotatebox{-90}{$\hookrightarrow$} &  &  &  & \rotatebox{-90}{$\hookrightarrow$} &  & \rotatebox{-90}{$\hookrightarrow$} &  & \rotatebox{-90}{$\hookrightarrow$} &  &  \\
        \cdots & \Omega_{n+1}^m &
        \xrightarrow[]{\partial_{q+1}^m} & \Omega_n^m &
        \xrightarrow[]{\partial_n^m} & \cdots &
        \xrightarrow[]{\partial_3^m} & \Omega_2^m & \xrightarrow[]{\partial_2^m} & \Omega_1^m & \xrightarrow[]{\partial_1^m} & \Omega_0^m & \xrightarrow[]{\partial_0^m} & \Omega_{-1}^m
    \end{array}\right.
\end{equation}
For the sake of simplicity, we use $\Omega_n^t$ to represent $\Omega_n(G_t)$. Suppose a subset of $\Omega_n^{s}$ whose boundary is in $\Omega_{n-1}^t$ as:
\begin{equation}
    \Omega_n^{t,s} \coloneqq \{ \alpha \in \Omega_n^{s} \ | \ \partial_n^{s}\alpha \in \Omega_{n-1}^{t}\}.
\end{equation}
The persistent $n$-boundary operator is denoted as $\eth_n^{t,s} : \Omega_n^{t,s} \to  \Omega_{n-1}^{t}$, and its corresponding adjoint operator is $(\eth_n^{t,s})^{\ast} : \Omega_{n-1}^{t}  \to  \Omega_n^{t,s}$. Therefore, the persistent $n$-th path Laplacian operator $\Delta_n^{t,s}: \Omega_n^t \to \Omega_n^t$ defined along the filtration is:
\begin{equation}
    \Delta_n^{t,s} = \eth_{n+1}^{t,s} \left(\eth_{n+1}^{t,s}\right)^\ast + \partial_n^{t^\ast} \partial_n^t.
\end{equation}
 Since $\Delta_n^{t,s}$ inherits the inner product from $\eth_{n+1}^{t,s}$, then the adjoint map $\left(\eth_{n+1}^{t,s}\right)^\ast$ is well defined. Intuitively, the matrix representation of $\Delta_n^{t,s}$ is
\begin{equation}
    L_n^{t,s} = B_{n+1}^{t,s} P^{-1} (B_{n+1}^{t,s})^T + (B_{n}^t)^T B_{n}^t,
\end{equation}
where $P^{-1}$ is the associated inner product matrix of $\Omega_{n+1}^{t,s}$ with arbitrary basis. Moreover, assume the dimension of $L_n^{t,s}$ is $N$, then the spectra of $L_n^{t,s}$ that are arranged in ascending order can be displayed as:
\[
\text{Spectra}(L_n^{t,s}) = \{(\lambda_1)_n^{t,s}, (\lambda_2)_n^{t,s}, \cdots, (\lambda_N)_n^{t,s}  \}.
\]
Note that the smallest non-harmonic spectra of $L_n^{t,s}$ is denoted as $(\tilde{\lambda}_2)_n^{t,s}$. We call the multiplicity of zero spectra of $L_q^{t,s}$ to be persistent $n$-th Betti number $\beta_n^{t,s}$ from $G_t$ to $G_s$.
\begin{equation}
    \beta_n^{t,s} = \text{nullity}(L_n^{t,s}) = \text{the number of zero eigenvalues (i.e., harmonic eigenvalues) of } L_n^{t,s}.
\end{equation}

{\bf Distanced-based filtration}  Specifically, suppose $G(w) = (V,E,w)$ is a weighted digraph, where $V$ is the set of the vertices and $E$ is the set of the directed edges. Assume $w$ is a weight function $w: E\to \mathbb{R}$. For example, if $V$ is in the Euclidean space, then a digraph $G(w)$ is a geometric digraph (a geometric digraph is a digraph in which the vertices are embedded as points in the Euclidean space, and the edges are embedded as non-crossing directed line segments). For any $(i,j) \in E$ where $i,j\in V$, we define $w(i,j) = \|i-j\|$, where $\|\cdot\|$ is a Euclidean metric. Hence, for every $\delta\in \mathbb{R}$, a digraph can be described as $G^{\delta} = (V, E^{\delta}) = (V,\{e\in E: w(e) \leq \delta\})$, and a filtration of digraphs can be described as $\{G^{\delta} \hookrightarrow G^{\delta^{\prime}}\}_{\delta \leq \delta^{\prime}}$.

Therefore, the persistent $n$-th path Laplacian matrix defined on the filtration is 
\begin{equation}
    L_n^{\delta,\delta^{\prime}} = B_{n+1}^{\delta,\delta^{\prime}} P^{-1} (B_{n+1}^{\delta,\delta^{\prime}})^T + (B_{n}^\delta)^T B_{n}^\delta,
\end{equation}
where its corresponding Betti numbers and spectra can be expressed as:
\begin{align}
    & \beta_n^{\delta,\delta^{\prime}} = \text{nullity}(L_n^{\delta,\delta^{\prime}}) = \text{the number of zero eigenvalues (i.e., harmonic eigenvalues) of } L_n^{\delta,\delta^{\prime}}.\\
    & \text{Spectra}(L_n^{\delta,\delta^{\prime}}) = \{(\lambda_1)_n^{\delta,\delta^{\prime}}, (\lambda_2)_n^{\delta,\delta^{\prime}}, \cdots, (\lambda_N)_n^{\delta,\delta^{\prime}} \}.
\end{align}
Notably, the Fiedler value (i.e., spectral gap) of $L_n^{\delta,\delta^{\prime}}$ is widely used in many other areas such as physics and geography, which is denoted as $\tilde{\lambda}_n^{\delta,\delta^{\prime}}$. As shown below, it is sensitive to both topological and geometric changes. 

 Moreover, it is worth to mention that isolated points (vertices) can be either included in the digraphs  (under the distance-based filtration)  or removed from the digraphs    (under the distanced-based filtration with removal of isolated points). 

\begin{figure}[ht!]
	\includegraphics[width=1.0\textwidth]{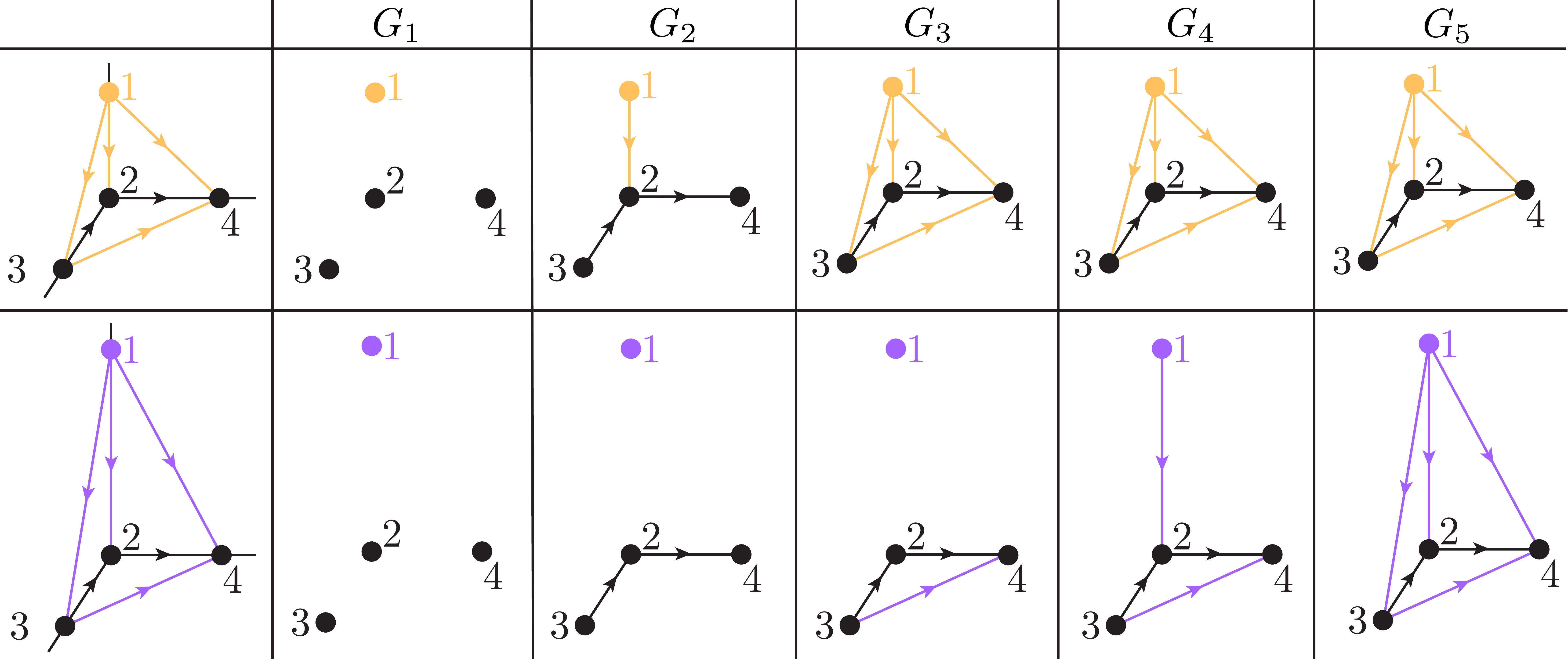}
	\centering
	\caption{Illustration of filtration on a tetrahedron. Here, $1,2,3,$ and $4$ represent four elementary $0$-paths $e_1, e_2, e_3,$ and $ e_4$. The top panel is a tetrahedron that has edge lengths $|e_{12}|= |e_{32}|= |e_{24}| = 1$ and $|e_{13}|= |e_{14}|=|e_{34}| = \sqrt{2}$. The bottom panel is a tetrahedron that has edge lengths $|e_{32}|= |e_{24}| = 1$, $|e_{34}|=\sqrt{2}$, $|e_{12}|=\sqrt{3}$, and $|e_{13}|=| e_{14}|=2$.}
	\label{fig:tetra persistent path laplacians}
\end{figure}

\begin{figure}[ht!]
	\includegraphics[width=1.0\textwidth]{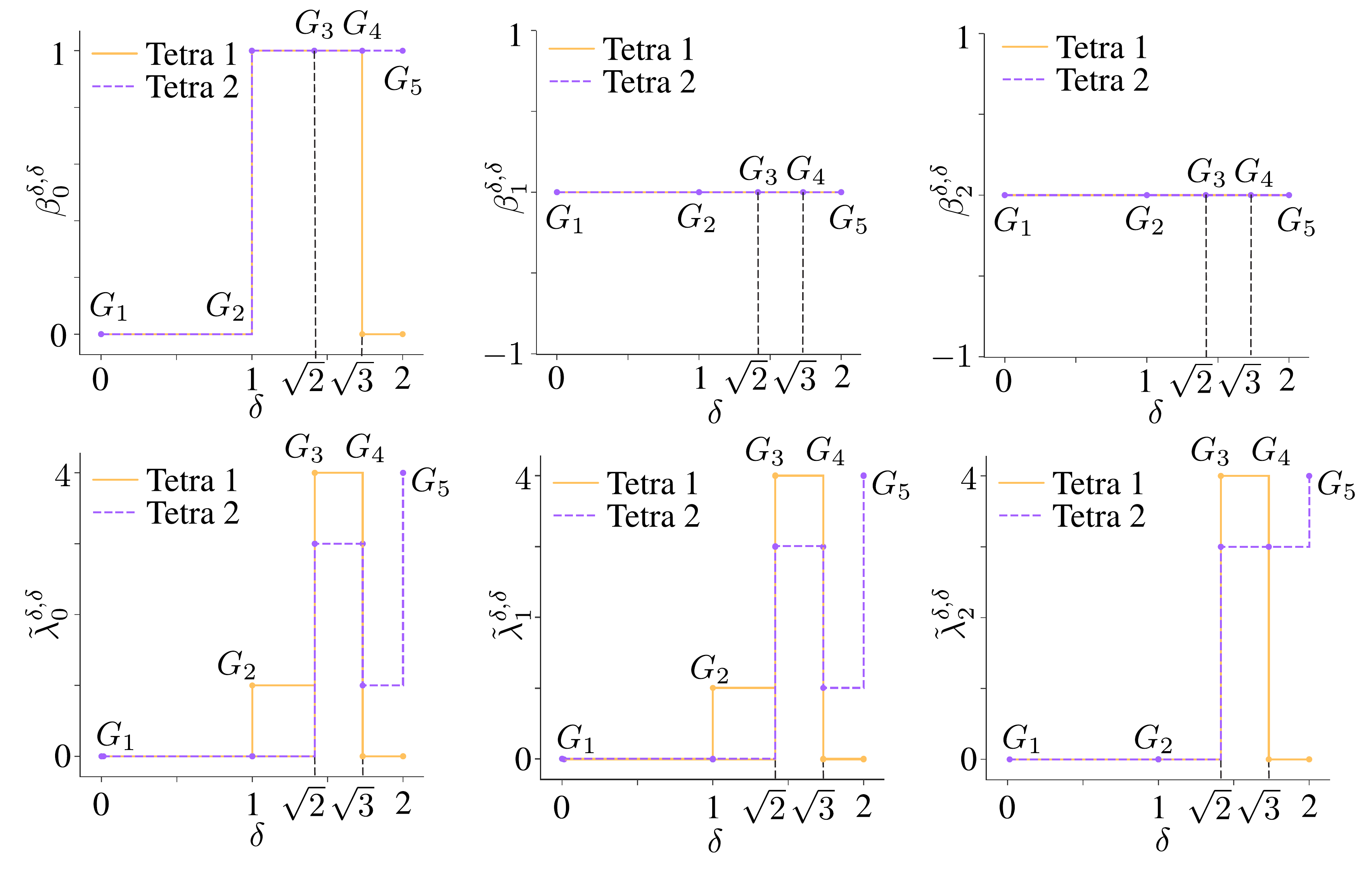}
	\centering
	\caption{Comparison of Betti numbers and non-harmonic spectra of $L_n^{\delta,\delta}$ when $n=0,1,$ and 2 on tetrahedrons Tetra 1 and Tetra 2. Note that 
since	$\beta_{1}^{\delta,\delta}=0$ and $\beta_{2}^{\delta,\delta}=0$ for Tetra 1 and Tetra 2, topological variants from persistent path homology cannot discriminate Tetra 1 and Tetra 2. However $\lambda_{1}^{\delta,\delta}$ and $\lambda_{2}^{\delta,\delta}$ show the differences between Tetra 1 and Tetra 2. 
	}
	\label{fig:tetra comparison}
\end{figure}

One can get both abstract information (revealed by Betti numbers) and geometric information (revealed by non-harmonic spectra) from digraphs along filtration. For instance, \autoref{fig:tetra persistent path laplacians} illustrates the filtration on two tetrahedrons. The top panel is a tetrahedron (Tetra 1) with edge lengths $|e_{12}|=| e_{32}|=| e_{24} |= 1$, and $|e_{13}|=| e_{14}|=| e_{34}| = \sqrt{2}$. The bottom panel is another tetrahedron (Tetra 2) with edge lengths $|e_{12}|=\sqrt{3}$, $|e_{32}|=| e_{24}| = 1$, and $|e_{13}|=| e_{14}|=2$, and $|e_{34} |= \sqrt{2}$. We say $G_1 = G^{0}, G_2 = G^{1}, G_3 = G^{\sqrt{2}}, G_4 = G^{\sqrt{3}},$ and $ G_5 = G^{\sqrt{5}}$. \autoref{fig:tetra comparison} shows the changes of $\beta_{n}^{\delta,\delta}$ and $\lambda_{n}^{\delta,\delta}$ of persistent $n$-th path Laplacian $L_n^{\delta,\delta}$ along filtration. It can be seen that by varying the filtration parameter $\delta$ from $0$ to $1$, the Betti 1 and Betti 2 are always 0. However, the smallest nonzero eigenvalue $\tilde{\lambda}_n^{\delta,\delta}$ of Tetra 1 and Tetra 2 have changes along filtration parameter $\delta$. Additionally, when $n=1, 2$, the $\tilde{\lambda}_n^{\delta,\delta}$ can distinguish Tetra 1 and Tetra 2, while $\beta_{n}^{\delta,\delta}$ cannot. This indicates that non-harmonic spectra of persistent path Laplacian can reveal more geometric information than the persistent Betti numbers in distinguishing similar topological structures.    Notably, we remove all the isolated points from each digraph for the simplicity of calculation.

\begin{figure}[ht!]
	\includegraphics[width=1.0\textwidth]{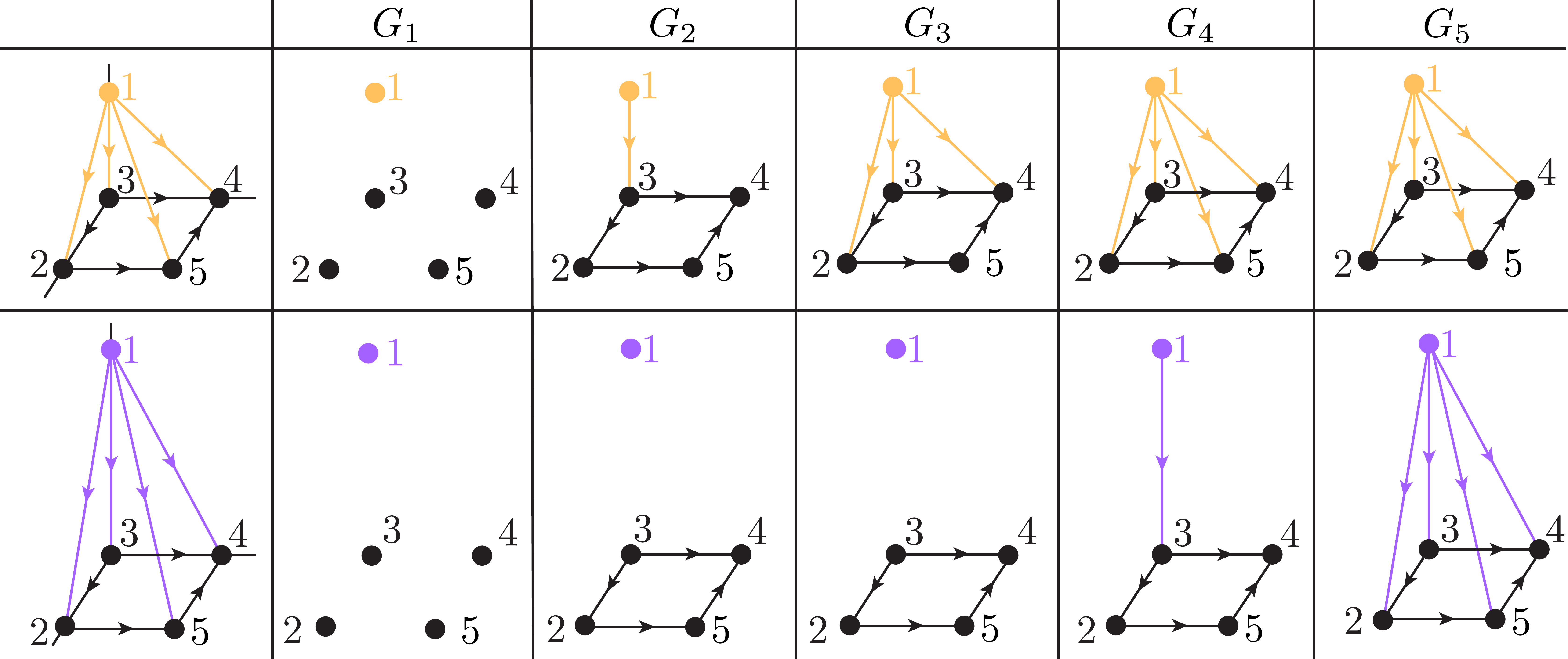}
	\centering
	\caption{Illustration of filtration on a pyramid. Here, $1,2,3,4,$ and $5$ represent five elementary $0$-paths $e_1, e_2, e_3,e_4,$ and $ e_5$. The top panel is a pyramid that has edge lengths $|e_{13}|=| e_{25}|=| e_{32}|=| e_{34}|=| e_{54}| = 1$, $|e_{12}|=|e_{14} |= \sqrt{2}$, and $|e_{15} |= \sqrt{3}$. The bottom panel is a pyramid that has edge lengths $|e_{25}|=|e_{32}|=|e_{34}|=|e_{54}| = 1$, $|e_{12}|=|e_{14}|= 2$, and $|e_{15} |= \sqrt{5}$.}
	\label{fig:pyramid persistent path laplacians}
\end{figure}

\begin{figure}[ht!]
	\includegraphics[width=1.0\textwidth]{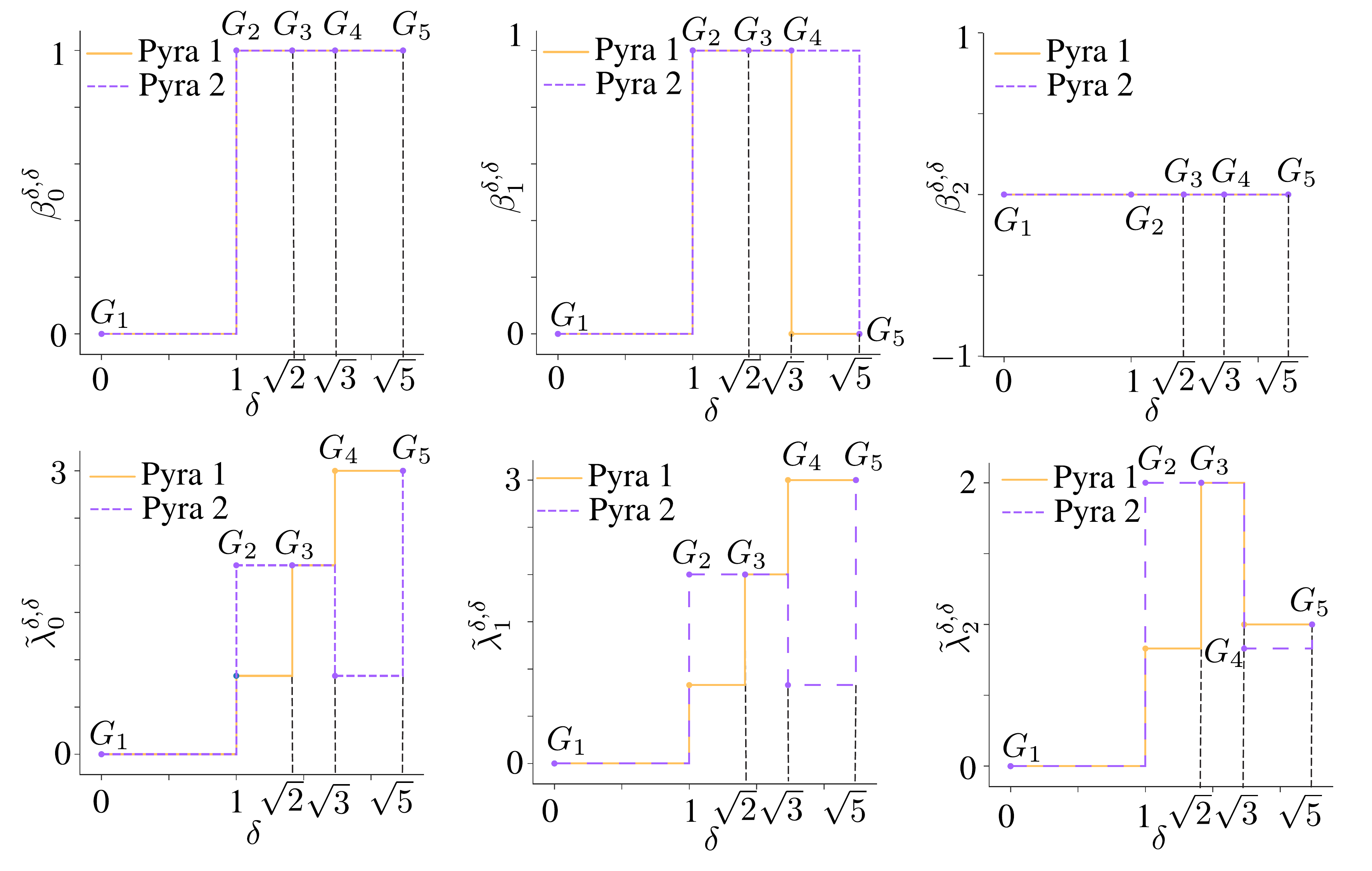}
	\centering
	\caption{Comparison of Betti number and non-harmonic spectra of $L_n^{\delta,\delta}$ when $n=0,1,$c and $2$ on pyramids Pyra 1 and Pyra 2.
Note that since	$\beta_{2}^{\delta,\delta}=0$, it cannot distinguish Pyra 1 and Pyra 2. But $\lambda_{2}^{\delta,\delta}$ can tell the difference.  
	}
	\label{fig:pyramid comparison}
\end{figure}

Moreover, a more complicated example is also illustrated in \autoref{fig:pyramid persistent path laplacians} to describe the filtration on two pyramids. The top panel is a pyramid (Pyra 1) with edge lengths $|e_{12}|=| e_{32}|=| e_{24}| = 1$, and $|e_{13}|=| e_{14}|=| e_{34}| = \sqrt{2}$. The bottom panel is a pyramid (Pyra 2) with edge lengths $|e_{12}|=\sqrt{3}$, $|e_{32},|=| e_{24}| = 1$, and $|e_{13}|=| e_{14}|=2$, and $|e_{34}| = \sqrt{2}$. We say $G_1 = G^{0}, G_2 = G^{1}, G_3 = G^{\sqrt{2}}, G_4 = G^{\sqrt{3}}, $ and $ G_5 = G^{\sqrt{5}}$. \autoref{fig:pyramid comparison} depicts the changes of $\beta_{n}^{\delta,\delta}$ and $\lambda_{n}^{\delta,\delta}$ of persistent $n$-th path Laplacian $L_n^{\delta,\delta}$ for objects Pyra 1 and Pyra 2 along filtration. For Pyra 1 and Pyra 2, when $n=0$ and $\delta=1$, their corresponding digraphs form, which result in $\beta_0^{1,1} = 1$ and $\beta_1^{1,1} = 1$ for both Pyra 1 and Pyra 2. When $\delta = \sqrt{3}$, we have  $\beta_1^{\sqrt{3},\sqrt{3}} = 0$ for Pyra 1 since the introducing of a new directed edges $e_{15}$. When $\delta = \sqrt{5}$, we have $\beta_1^{\sqrt{5},\sqrt{5}} = 0$ for Pyra 2 since the introducing of a new directed edges $e_{15}$ kills the $1$-cycle formed by $e_{25}, e_{32}, e_{34},$ and $ e_{54}$. Furthermore, although Pyra 1 and Pyra 2 do not have exactly the same geometric structure, their share the same $\beta_2^{\delta,\delta}$ value from $\delta = 0$ to $\delta = \sqrt{5}$.  However, Pyra 1 and Pyra 2 can be distinguished by the $\tilde{\lambda}_2^{\delta,\delta}$ along filtration. Therefore, we can see that similar to the PSG, one can use the non-harmonic spectra from the persistent path laplacian to reveal the intrinsic geometric information of a givens point-cloud dataset by varying the filtration parameters. In addition, the detailed calculations of $L_n^{\delta,\delta}$ can be found in the Appendix. 

\section{Application}
In this section, we apply the persistent path Laplacian to the analysis of the curcurbit[n]urils system. Cucurbiturils are macrocyclic molecules, which are made of glycoluril(=C$_6$H$_2$N$_4$O$_2$=) monomers linked by methylene bridges (-CH$_2$-). CB$n$ is commonly used  as an abbreviation of Cucurbiturils. Here, $n$ is the number of glycoluril units. In this work, we consider CB7 as an example. 
The  molecular formulas of CB7 is C$_{42}$H$_{14}$N$_{28}$O$_{14}$. 
The molecular structure of CB7 is obtained from the Supporting Information of Ref. \cite{gao2015binding}. 

\begin{figure}[H]
    \centering
    \includegraphics[width=1\textwidth]{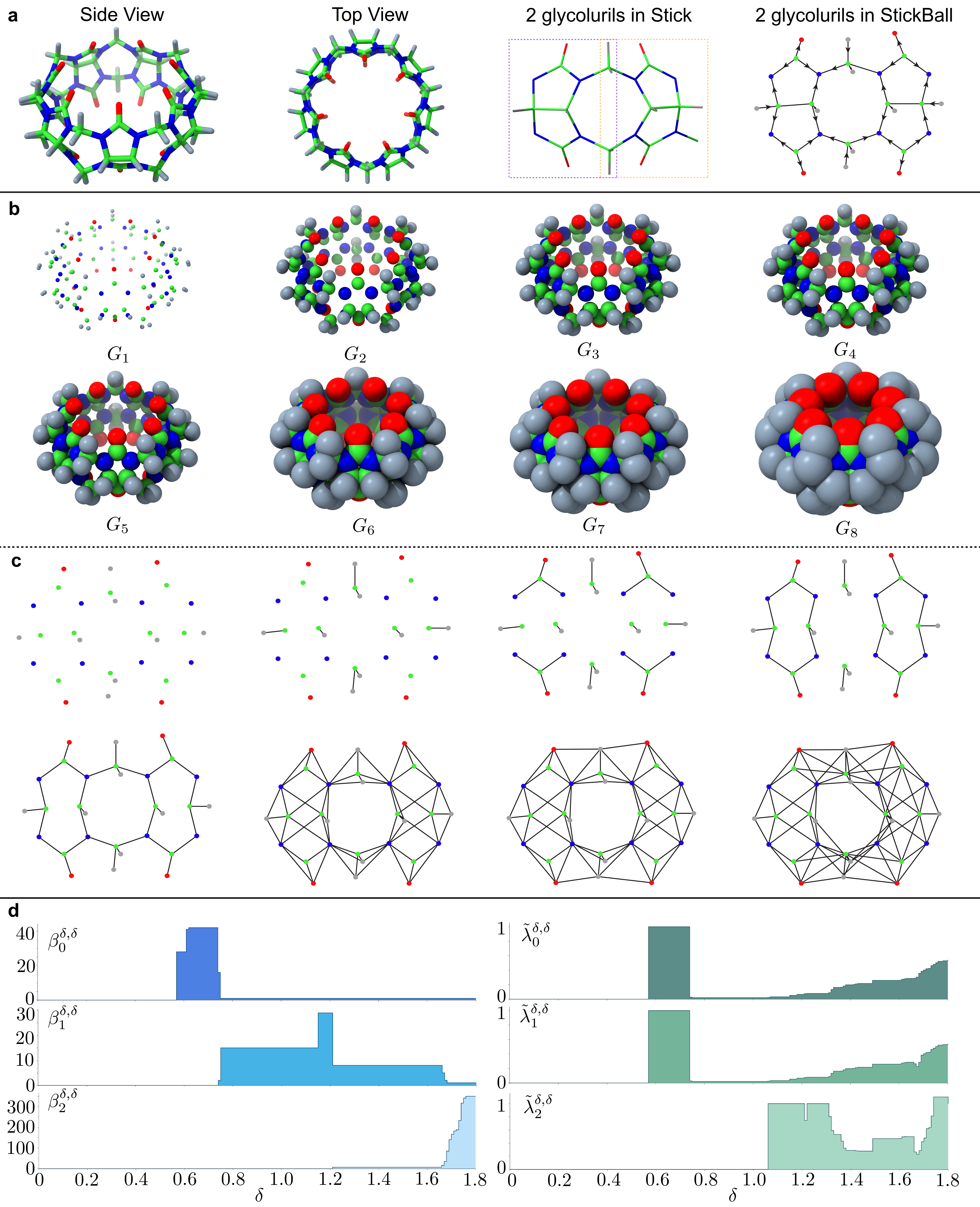}
    \caption{{\bf a} The 3D structures of CB7, 2 glycolurils, and path direction assignment. Here, from left to right, the side view of CB7, top view of CB7,  the structure of two glycoluril units (=C$_{10}$H$_4$N$_8$O$_4$=), and electronegativity-based  path direction assignment are depicted as well. {\bf b} Illustration of filtration-induced geometries $G_i (i=1,2,...,8)$ of CB7. Eight digraphs $G_1 = G_0^{0.200}, G_2 = G_0^{0.565}, G_3 = G_0^{0.710}, G_4 = G_0^{0.745}, G_5 = G_0^{0.800}, G_6 = G_0^{1.210}, G_7 = G_0^{1.315}, G_8=G_0^{1.800}$ are constructed under filtration parameter $\delta$. {\bf c} Illustration of filtration-induced  path complexes within two glycoluril units.  Path directions can be inferred from their colors as shown in the last chart of {\bf a}.  
		{\bf d} Betti numbers $\beta_n^{\delta,\delta}$ and non-harmonic spectra $\tilde{\lambda}_n^{\delta,\delta}$ of persistent path Laplacians ($L_n^{\delta,\delta}$ when $n=0,1,$ and $2$) for CB7. }
    \label{fig:CB7}
\end{figure}

\autoref{fig:CB7} illustrates how PPL is employed for a molecular system to extract its rich topological and geometric information. The first two charts of \autoref{fig:CB7}{\bf a} describe the three-dimensional (3D) top view and side view of CB7. The green, blue, red, and gray colors represent  C, N, O, and H atoms, respectively. The third chart of \autoref{fig:CB7}{\bf a} is a basic ``Octagon-pentagon'' unit that consists of two glycolurils. It can be seen that $7$ glycolurils exist in CB7. The last chart of \autoref{fig:CB7}{\bf a} demonstrates the path direction assignment to   pairs of atoms based on atomic electronegativity. The periodic table of electronegativity is given by the Pauling scale \cite{pauling1932nature}, in which the  electronegativities of C, N, O, and H are 2.55, 3.04, 3.44, and 2.20, respectively. Then, we set the directions of edges following the order ``$\text{H} \to \text{C} \to \text{N} \to \text{O}$". 

\autoref{fig:CB7}{\bf b} depicts the distance-based filtration of CB7. Here, structures $G_i (i=1,2,..., 8)$ were obtained at the filtration radii of 0.200, 0.565, 0.710, 0.745, 0.800, 1.210, 1.315, and  1.800 {\AA}, respectively. In our digraph notation, we denote these structures as $G_1 = G_0^{0.200}, G_2 = G_0^{0.565},G_3 = G_0^{0.710}, G_4 = G_0^{0.745}$, $G_5 = G_0^{0.800}, G_6 = G_0^{1.210}, G_7 = G_0^{1.315}, $ and $ G_8=G_0^{1.800}$. Note that, in the present formulation,  all of the isolated points were removed from these digraphs. 

\autoref{fig:CB7}{\bf c} illustrates the filtration-induced path complexes in the aforementioned $G_i  (i=1,2,..., 8)$. To clearly show the topological and geometric changes, only the path complexes in one ``Octagon-pentagon'' unit (or two glycolurils) are considered and depicted for each structure. For simplicity, only edges are presented. However, their path directions can be easily assigned based on their color map as shown in the last chart of  \autoref{fig:CB7}{\bf a}.  

\autoref{fig:CB7}{\bf d} depicts the PPL spectra of CB7. We can see that at the initial state ($G_1$) when $\delta = \SI{0.200}{\angstrom}$ ), total $98$ atoms  are isolated from one another. When radius $\delta = \SI{0.565}{\angstrom}$ ($G_2$),  C atoms on each pentagon are connected with their H atom neighborhoods. Therefore, four isolated components are formed in each glycoluril, which makes $\beta_{0}^{\delta,\delta} = 4 \times 7 = 28$. 
At $G_3$ ($r= \SI{0.710}{\angstrom}$), C atoms on each pentagon are connected with their N and O neighborhoods. At this stage, two more connected components are involved in one glycoluri structure, which makes $\beta_{0}^{\delta,\delta} = 6 \times 7 = 42$. Only one connected structure can be formed if all of the atoms get connected with their neighborhood atoms. Therefore, $\beta_{0}^{\delta,\delta} = 1$ (see $G_5$ - $G_8$). Notably, the $\beta_{2}^{\delta,\delta}$ and $\tilde{\lambda}_2^{\delta,\delta}$ provide rich topological and geometric information when the filtration parameter $\delta$ increases.

This example shows that PPL can decode topological persistence and the shape evolution of a given molecular system with chemical- or biological-based directional assignment. Specifically,   $\tilde{\lambda}_0^{\delta,\delta}$ can still offer geometric information when $\beta_{0}^{\delta,\delta}$ does not changes for large radii. Therefore, PPL keeps revealing homotopic shape evolution when the topological invariant from persistent path homology does not change.  

Additionally, unlike persistent Laplacian, high-order PPL operators provide rich topological information. For instance, when the filtration parameter $\delta$ increases to $1.68$, $\beta_{2}^{\delta,\delta}$ from PPL dramatically goes up. Whereas,  in persistent Laplacian, the value of Betti 2 is quite limited since the CB7 system can barely form $2$-cycles at a similar filtration parameter using either Rips complex or alpha complex. This trait endows PPL with a better ability to characterize the geometry and topology of an object at large scales.   
 

\section{Conclusion}
 
Path homology, a rich mathematical concept introduced by Grigor’yan, Lin, Muranov, and Yau, has stimulated a variety of new developments in pure and applied mathematics, including much attention from the topological data analysis (TDA) community. Unlike original homology or persistent homology, path homology enables the treatment of directed graphs and networks. Persistent path homology bridges path homology with multiscale analysis,  making it a powerful tool for practical applications. Nonetheless, these formulations are insensitive to homotopic shape evolution during filtration. 

Topological Laplacians, including Hodge Laplacian, graph Laplacian, sheaf Laplacian, and Dirac Laplacian,  are versatile mathematical tools that not only preserve all topological invariants but also describe geometric shapes. This work introduces a new topological Laplacian, namely persistent path Laplacian,  as a new mathematical tool for the multi-scale analysis of  directed graphs and networks. For a given data, the proposed persistent path Laplacian fully recovers the topological persistence of persistent homology in its harmonic spectra and meanwhile, captures homotopic shape evolution of the data during filtration in its non-harmonic spectra. 

\section*{Supporting information} 
Supporting information is given on the detailed boundary matrices, Laplacian matrices, Betti numbers and eigenvalues for the digraphs shown in    
\autoref{fig:pyramid persistent path laplacians}.

\section*{Acknowledgment}
This work was supported in part by NIH grants  R01GM126189 and  R01AI164266, NSF grants DMS-2052983,  DMS-1761320, and IIS-1900473,  NASA grant 80NSSC21M0023,  Michigan Economic Development Corporation, MSU Foundation,  Bristol-Myers Squibb 65109, and Pfizer.

\section*{Appendix}

In Tables 5-19, we present the detailed matrix constructions, Betti numbers, and spectra for various digraphs shown in Figure 5 top and bottom panels. 

\begin{table}[H]
    \centering
    \setlength\tabcolsep{6pt}
    \captionsetup{margin=0.9cm}
    \caption{Matrix construction of graph $G_1$ (with isolated points included) in the top panel of Figure 5.}
    \begin{tabular}{c|ccc}
    \hline
    $n$  & $n=0$ & $n=1$ & $n=2$ \\ \hline 
    $\Omega_n$ & $\text{span}\{e_1,e_2,e_3,e_4, e_5\}$ & $\{0\}$ & $\{0\}$ \\ \hline \\
    $B_{n+1}$  & $5\times 0$ empty matrix
        & / &  /    \\  \\
    $L_n$  & $5\times 5$ zero matrix  & / &    /  \\  \\
    $\beta_{n}$                        & 5               & / & /    \\  \\
    $\text{Spectra}(L_{n})$    & $\{0,0,0,0,0\}$   & / & /  \\ \hline
    \end{tabular}
    \\
    \label{tab:Figure 5 G1}
\end{table}

\begin{table}[H]
    \centering
    \setlength\tabcolsep{6pt}
    \captionsetup{margin=0.9cm}
    \caption{Matrix construction of graph $G_1$ (without isolated points) in the top panel of Figure 5.}
    \begin{tabular}{c|ccc}
    \hline
    $n$  & $n=0$ & $n=1$ & $n=2$ \\ \hline 
    $\Omega_n$ & $\{0\}$ & $\{0\}$ & $\{0\}$ \\ \hline \\
    $B_{n+1}$  & /
        & / &  /    \\  \\
    $L_n$  & /  & / &    /  \\  \\
    $\beta_{n}$                        & /               & / & /    \\  \\
    $\text{Spectra}(L_{n})$    & /   & / & /  \\ \hline
    \end{tabular}
    \\
    \label{tab:Figure 5 G1}
\end{table}

\begin{table}[H]
    \centering
    \setlength\tabcolsep{6pt}
    \captionsetup{margin=0.9cm}
    \caption{Matrix construction of graph $G_2$ in the top panel of Figure 5.}
    \begin{tabular}{c|ccc}
    \hline
    $n$  & $n=0$ & $n=1$ & $n=2$ \\ \hline 
    $\Omega_n$ & $\text{span}\{e_1,e_2,e_3,e_4,e_5\}$ & $\text{span}\{e_{13}, e_{25}, e_{32}, e_{34},e_{45}\}$ & $\{0\}$ \\ \hline \\
    $B_{n+1}$  & $\begin{array}{@{}r@{}c@{}c@{}c@{}c@{}c@{}l@{}}
            & e_{13} & e_{25} & e_{32} & e_{34} & e_{45}  \\
           \left.\begin{array}{c}
            e_1 \\
            e_2 \\
            e_3 \\
            e_4 \\
            e_5
            \end{array}\right(
            & \begin{array}{c}-1 \\  0 \\  1  \\  0  \\  0  \end{array}
            & \begin{array}{c} 0 \\ -1 \\  0  \\  0  \\  1  \end{array}
            & \begin{array}{c} 0 \\  1 \\ -1  \\  0  \\  0  \end{array}
            & \begin{array}{c} 0 \\  0 \\ -1  \\  1  \\  0  \end{array}
            & \begin{array}{c} 0 \\  0 \\  0  \\  1  \\ -1  \end{array}
            & \left)\begin{array}{c} \\ \\  \\ \\ \\
            \end{array}\right.
        \end{array}$
        & $5\times 0$ empty matrix &  $\left(\begin{array}{ccccc}
         /
    \end{array}\right)$   \\  \\
    $L_n$  & $\left(\begin{array}{ccccc}
         1  &   0  &  -1  &   0  &   0  \\
         0  &   2  &  -1  &   0  &  -1  \\
        -1  &  -1  &   3  &  -1  &   0  \\
         0  &   0  &  -1  &   2  &  -1  \\
         0  &  -1  &   0  &  -1  &   2   
    \end{array}\right)$  & $\left(\begin{array}{ccccc}
         2  &   0  &  -1  &  -1  &   0  \\
         0  &   2  &  -1  &   0  &  -1  \\
        -1  &  -1  &   2  &   1  &   0  \\
        -1  &   0  &   1  &   2  &   1  \\
         0  &  -1  &   0  &   1  &   2   
    \end{array}\right)$ &    $\left(\begin{array}{ccccc}
         / 
    \end{array}\right)$  \\  \\
    $\beta_{n}$                        & 1               & 1 & 0    \\  \\
    $\text{Spectra}(L_{n})$    & $\{0,0.8299,2,2.6889,4.4812\}$   & $\{0,0.8299,2,2.6889,4.4812\}$ & / \\ \hline
    \end{tabular}
    \\
    \label{tab:Figure 5 G2}
\end{table}

\begin{table}[H]\footnotesize
    \centering
    \setlength\tabcolsep{6pt}
    \captionsetup{margin=0.9cm}
    \caption{Matrix construction of graph $G_3$ in the top panel of Figure 5.}
    \begin{tabular}{c|ccc}
    \hline
    $n$  & $n=0$ & $n=1$ & $n=2$ \\ \hline 
    $\Omega_n$ & $\text{span}\{e_1,e_2,e_3,e_4,e_5\}$ & $\text{span}\{e_{12}, e_{13}, e_{14}, e_{25}, e_{32}, e_{34}, e_{54}\}$ & $\text{span}\{e_{132},e_{134}\}$ \\ \hline \\
    $B_{n+1}$  & $\begin{array}{@{}r@{}c@{}c@{}c@{}c@{}c@{}c@{}c@{}l@{}}
            & e_{12} & e_{13} & e_{14} & e_{25} & e_{32} & e_{34} & e_{54}  \\
           \left.\begin{array}{c}
            e_1 \\
            e_2 \\
            e_3 \\
            e_4 \\
            e_5
            \end{array}\right(
            & \begin{array}{c} -1 \\  1  \\  0  \\  0  \\ 0  \end{array}
            & \begin{array}{c} -1 \\  0  \\  1  \\  0  \\ 0  \end{array}
            & \begin{array}{c} -1 \\  0  \\  0  \\  1  \\ 0  \end{array}
            & \begin{array}{c}  0 \\ -1  \\  0  \\  0  \\ 1  \end{array}
            & \begin{array}{c}  0 \\  1  \\ -1  \\  0  \\ 0  \end{array}
            & \begin{array}{c}  0 \\  0  \\ -1  \\  1  \\ 0  \end{array}
            & \begin{array}{c}  0 \\  0  \\  0  \\  1  \\-1  \end{array}
            & \left)\begin{array}{c} \\ \\  \\ \\ \\ 
            \end{array}\right.
        \end{array}$
        & $\begin{array}{@{}r@{}c@{}c@{}c@{}c@{}c@{}c@{}c@{}l@{}}
            & e_{132} & e_{134}  \\
           \left.\begin{array}{c}
            e_{12} \\
            e_{13} \\
            e_{14} \\
            e_{25} \\
            e_{32} \\
            e_{34} \\
            e_{54}
            \end{array}\right(
            & \begin{array}{c} -1 \\  1  \\  0  \\  0  \\ 1 \\ 0 \\ 0  \end{array}
            & \begin{array}{c}  0 \\  1  \\ -1  \\  0  \\ 0 \\ 1 \\ 0  \end{array}
            & \left)\begin{array}{c} \\ \\  \\ \\ \\ \\ \\
            \end{array}\right.
        \end{array}$ &  $2\times 0$ empty matrix  \\  \\
    $L_n$  & $\left(\begin{array}{ccccc}
         3  &  -1  &  -1  &  -1  &  0  \\
        -1  &   3  &  -1  &   0  & -1  \\
        -1  &  -1  &   3  &  -1  &  0  \\
        -1  &   0  &  -1  &   3  & -1  \\
         0  &  -1  &   0  &  -1  & 2  
    \end{array}\right)$  & $\left(\begin{array}{cccccccccc}
         3  &   0  &   1  &  -1  &  0  &  0  &  0 \\
         0  &   4  &   0  &   0  &  0  &  0  &  0 \\
         1  &   0  &   3  &   0  &  0  &  0  &  0 \\
        -1  &   0  &   0  &   2  & -1  &  0  & -1 \\
         0  &   0  &   0  &  -1  &  3  &  1  &  0 \\
         0  &   0  &   0  &   0  &  1  &  3  &  1 \\
         0  &   0  &   1  &  -1  &  0  &  1  &  2 
    \end{array}\right)$ &    $\left(\begin{array}{ccccc}
         3 & 1 \\
         1 & 3
    \end{array}\right)$  \\  \\
    $\beta_{n}$                        & 1               & 1 & 0    \\  \\
    $\text{Spectra}(L_{n})$    & $\{0,2,3,4,5\}$   & $\{0,2,2,3,4,4,5\}$ & $\{2,4\}$ \\ \hline
    \end{tabular}
    \\
    \label{tab:Figure 5 G3}
\end{table}

\begin{landscape}

\begin{table}[H]\footnotesize
    \centering
    \setlength\tabcolsep{6pt}
    \captionsetup{margin=0.9cm}
    \caption{Matrix construction of graph $G_4$ in the top panel of Figure 5.}
    \begin{tabular}{c|ccc}
    \hline
    $n$  & $n=0$ & $n=1$ & $n=2$ \\ \hline 
    $\Omega_n$ & $\text{span}\{e_1,e_2,e_3,e_4,e_5\}$ & $\text{span}\{e_{12}, e_{13}, e_{14}, e_{15}, e_{25}, e_{32}, e_{34}, e_{54}\}$ & $\text{span}\{e_{125},e_{132},e_{134},e_{154}\}$ \\ \hline \\
    $B_{n+1}$  & $\begin{array}{@{}r@{}c@{}c@{}c@{}c@{}c@{}c@{}c@{}c@{}l@{}}
            & e_{12} & e_{13} & e_{14} & e_{15} & e_{25} & e_{32} & e_{34} & e_{54}  \\
           \left.\begin{array}{c}
            e_1 \\
            e_2 \\
            e_3 \\
            e_4 \\
            e_5
            \end{array}\right(
            & \begin{array}{c} -1 \\  1  \\  0  \\  0  \\ 0  \end{array}
            & \begin{array}{c} -1 \\  0  \\  1  \\  0  \\ 0  \end{array}
            & \begin{array}{c} -1 \\  0  \\  0  \\  1  \\ 0  \end{array}
            & \begin{array}{c} -1 \\  0  \\  0  \\  0  \\ 1  \end{array}
            & \begin{array}{c}  0 \\ -1  \\  0  \\  0  \\ 1  \end{array}
            & \begin{array}{c}  0 \\  1  \\ -1  \\  0  \\ 0  \end{array}
            & \begin{array}{c}  0 \\  0  \\ -1  \\  1  \\ 0  \end{array}
            & \begin{array}{c}  0 \\  0  \\  0  \\  1  \\-1  \end{array}
            & \left)\begin{array}{c} \\ \\  \\ \\ \\ 
            \end{array}\right.
        \end{array}$
        & $\begin{array}{@{}r@{}c@{}c@{}c@{}c@{}c@{}c@{}c@{}l@{}}
            & e_{125} & e_{132} & e_{134} & e_{154} \\
           \left.\begin{array}{c}
            e_{12} \\
            e_{13} \\
            e_{14} \\
            e_{15} \\
            e_{25} \\
            e_{32} \\
            e_{34} \\
            e_{54}
            \end{array}\right(
            & \begin{array}{c}  1 \\  0  \\  0  \\ -1  \\ 1 \\ 0 \\ 0  \\ 0   \end{array}
            & \begin{array}{c} -1 \\  1  \\  0  \\  0  \\ 0 \\ 1 \\ 0  \\ 0   \end{array}
            & \begin{array}{c}  0 \\  1  \\ -1  \\  0  \\ 0 \\ 0 \\ 1  \\ 0   \end{array}
            & \begin{array}{c}  0 \\  0  \\ -1  \\  1  \\ 0 \\ 0 \\ 0  \\ 1   \end{array}
            & \left)\begin{array}{c} \\ \\  \\ \\ \\ \\ \\ \\
            \end{array}\right.
        \end{array}$ &  $4\times 0$ empty matrix  \\  \\
    $L_n$  & $\left(\begin{array}{ccccc}
         4  &  -1  &  -1  &  -1  &  -1  \\
        -1  &   3  &  -1  &   0  & -1  \\
        -1  &  -1  &   3  &  -1  &  0  \\
        -1  &   0  &  -1  &   3  & -1  \\
        -1  &  -1  &   0  &  -1  &  3  
    \end{array}\right)$  & $\left(\begin{array}{cccccccccc}
         4  &   0  &   1  &   0  &  0  &  0  &  0  &  0 \\
         0  &   4  &   0  &   1  &  0  &  0  &  0  &  0 \\
         1  &   0  &   4  &   0  &  0  &  0  &  0  &  0 \\
         0  &   1  &   0  &   4  &  0  &  0  &  0  &  0 \\
         0  &   0  &   0  &   0  &  3  & -1  &  0  & -1 \\
         0  &   0  &   0  &   0  & -1  &  3  &  1  &  0 \\
         0  &   0  &   0  &   0  &  0  &  1  &  3  &  1 \\
         0  &   0  &   0  &   0  & -1  &  0  &  1  &  3 \\
    \end{array}\right)$ &    $\left(\begin{array}{ccccc}
         3 & -1 & 0 & -1 \\
        -1 &  3 & 1 &  0 \\
         0 &  1 & 3 &  1 \\
        -1 &  0 & 1 &  3
    \end{array}\right)$  \\  \\
    $\beta_{n}$                        & 1               & 1 & 0    \\  \\
    $\text{Spectra}(L_{n})$    & $\{0,3,3,5,5\}$   & $\{1,3,3,3,3,5,5,5\}$ & $\{1,3,3,5\}$ \\ \hline
    \end{tabular}
    \\
    \label{tab:Figure 5 G4}
\end{table}
\end{landscape}

\begin{landscape}
\begin{table}[H]\footnotesize
    \centering
    \setlength\tabcolsep{6pt}
    \captionsetup{margin=0.9cm}
    \caption{Matrix construction of graph $G_5$ in the top panel of Figure 5.}
    \begin{tabular}{c|ccc}
    \hline
    $n$  & $n=0$ & $n=1$ & $n=2$ \\ \hline 
    $\Omega_n$ & $\text{span}\{e_1,e_2,e_3,e_4,e_5\}$ & $\text{span}\{e_{12}, e_{13}, e_{14}, e_{15}, e_{25}, e_{32}, e_{34}, e_{54}\}$ & $\text{span}\{e_{125},e_{132},e_{134},e_{154}\}$ \\ \hline \\
    $B_{n+1}$  & $\begin{array}{@{}r@{}c@{}c@{}c@{}c@{}c@{}c@{}c@{}c@{}l@{}}
            & e_{12} & e_{13} & e_{14} & e_{15} & e_{25} & e_{32} & e_{34} & e_{54}  \\
           \left.\begin{array}{c}
            e_1 \\
            e_2 \\
            e_3 \\
            e_4 \\
            e_5
            \end{array}\right(
            & \begin{array}{c} -1 \\  1  \\  0  \\  0  \\ 0  \end{array}
            & \begin{array}{c} -1 \\  0  \\  1  \\  0  \\ 0  \end{array}
            & \begin{array}{c} -1 \\  0  \\  0  \\  1  \\ 0  \end{array}
            & \begin{array}{c} -1 \\  0  \\  0  \\  0  \\ 1  \end{array}
            & \begin{array}{c}  0 \\ -1  \\  0  \\  0  \\ 1  \end{array}
            & \begin{array}{c}  0 \\  1  \\ -1  \\  0  \\ 0  \end{array}
            & \begin{array}{c}  0 \\  0  \\ -1  \\  1  \\ 0  \end{array}
            & \begin{array}{c}  0 \\  0  \\  0  \\  1  \\-1  \end{array}
            & \left)\begin{array}{c} \\ \\  \\ \\ \\ 
            \end{array}\right.
        \end{array}$
        & $\begin{array}{@{}r@{}c@{}c@{}c@{}c@{}c@{}c@{}c@{}l@{}}
            & e_{125} & e_{132} & e_{134} & e_{154} \\
           \left.\begin{array}{c}
            e_{12} \\
            e_{13} \\
            e_{14} \\
            e_{15} \\
            e_{25} \\
            e_{32} \\
            e_{34} \\
            e_{54}
            \end{array}\right(
            & \begin{array}{c}  1 \\  0  \\  0  \\ -1  \\ 1 \\ 0 \\ 0  \\ 0   \end{array}
            & \begin{array}{c} -1 \\  1  \\  0  \\  0  \\ 0 \\ 1 \\ 0  \\ 0   \end{array}
            & \begin{array}{c}  0 \\  1  \\ -1  \\  0  \\ 0 \\ 0 \\ 1  \\ 0   \end{array}
            & \begin{array}{c}  0 \\  0  \\ -1  \\  1  \\ 0 \\ 0 \\ 0  \\ 1   \end{array}
            & \left)\begin{array}{c} \\ \\  \\ \\ \\ \\ \\ \\
            \end{array}\right.
        \end{array}$ &  $4\times 0$ empty matrix  \\  \\
    $L_n$  & $\left(\begin{array}{ccccc}
         4  &  -1  &  -1  &  -1  &  -1  \\
        -1  &   3  &  -1  &   0  & -1  \\
        -1  &  -1  &   3  &  -1  &  0  \\
        -1  &   0  &  -1  &   3  & -1  \\
        -1  &  -1  &   0  &  -1  &  3  
    \end{array}\right)$  & $\left(\begin{array}{cccccccccc}
         4  &   0  &   1  &   0  &  0  &  0  &  0  &  0 \\
         0  &   4  &   0  &   1  &  0  &  0  &  0  &  0 \\
         1  &   0  &   4  &   0  &  0  &  0  &  0  &  0 \\
         0  &   1  &   0  &   4  &  0  &  0  &  0  &  0 \\
         0  &   0  &   0  &   0  &  3  & -1  &  0  & -1 \\
         0  &   0  &   0  &   0  & -1  &  3  &  1  &  0 \\
         0  &   0  &   0  &   0  &  0  &  1  &  3  &  1 \\
         0  &   0  &   0  &   0  & -1  &  0  &  1  &  3 \\
    \end{array}\right)$ &    $\left(\begin{array}{ccccc}
         3 & -1 & 0 & -1 \\
        -1 &  3 & 1 &  0 \\
         0 &  1 & 3 &  1 \\
        -1 &  0 & 1 &  3
    \end{array}\right)$  \\  \\
    $\beta_{n}$                        & 1               & 0 & 0    \\  \\
    $\text{Spectra}(L_{n})$    & $\{0,3,3,5,5\}$   & $\{1,3,3,3,3,5,5,5\}$ & $\{1,3,3,5\}$ \\ \hline
    \end{tabular}
    \\
    \label{tab:Figure 5 G5}
\end{table}
\end{landscape}

\begin{table}[H]
    \centering
    \setlength\tabcolsep{6pt}
    \captionsetup{margin=0.9cm}
    \caption{Matrix construction of graph $G_1$ (with isolated points included) in the bottom panel of Figure 5.}
    \begin{tabular}{c|ccc}
    \hline
    $n$  & $n=0$ & $n=1$ & $n=2$ \\ \hline 
    $\Omega_n$ & $\text{span}\{e_1,e_2,e_3,e_4, e_5\}$ & / & / \\ \hline \\
    $B_{n+1}$  & $5\times 0$ empty matrix
        & / &  /    \\  \\
    $L_n$  & $5\times 5$ zero matrix  & / &    /  \\  \\
    $\beta_{n}$                        & 5               & / & /    \\  \\
    $\text{Spectra}(L_{n})$    & $\{0,0,0,0,0\}$   & / & /  \\ \hline
    \end{tabular}
    \\
    \label{tab:Figure 5 bottom G1}
\end{table}

\begin{table}[H]
    \centering
    \setlength\tabcolsep{6pt}
    \captionsetup{margin=0.9cm}
    \caption{Matrix construction of graph $G_1$ (without isolated points) in the bottom panel of Figure 5.}
    \begin{tabular}{c|ccc}
    \hline
    $n$  & $n=0$ & $n=1$ & $n=2$ \\ \hline 
    $\Omega_n$ & $\{0\}$ & $\{0\}$ & $\{0\}$ \\ \hline \\
    $B_{n+1}$  & /
        & / &  /    \\  \\
    $L_n$  & /  & / &    /  \\  \\
    $\beta_{n}$                        & /               & / & /    \\  \\
    $\text{Spectra}(L_{n})$    & /   & / & /  \\ \hline
    \end{tabular}
    \\
    \label{tab:Figure 5 bottom G1}
\end{table}

\begin{table}[H]
    \centering
    \setlength\tabcolsep{6pt}
    \captionsetup{margin=0.9cm}
    \caption{Matrix construction of graph $G_2$ (with isolated points included) in the bottom panel of Figure 5.}
    \begin{tabular}{c|ccc}
    \hline
    $n$  & $n=0$ & $n=1$ & $n=2$ \\ \hline 
    $\Omega_n$ & $\text{span}\{e_1,e_2,e_3,e_4,e_5\}$ & $\text{span}\{e_{25}, e_{32}, e_{34}, e_{54}\}$ & $\{0\}$ \\ \hline \\
    $B_{n+1}$  & $\begin{array}{@{}r@{}c@{}c@{}c@{}c@{}c@{}l@{}}
            & e_{25} & e_{32} & e_{34} & e_{54}  \\
           \left.\begin{array}{c}
            e_1 \\
            e_2 \\
            e_3 \\
            e_4 \\
            e_5
            \end{array}\right(
            & \begin{array}{c}  0 \\ -1  \\  0  \\  0  \\  1  \end{array}
            & \begin{array}{c}  0 \\  1  \\ -1  \\  0  \\  0  \end{array}
            & \begin{array}{c}  0 \\  0  \\ -1  \\  1  \\  0  \end{array}
            & \begin{array}{c}  0 \\  0  \\  0  \\  1  \\ -1  \end{array}
            & \left)\begin{array}{c} \\ \\  \\ \\ \\
            \end{array}\right.
        \end{array}$
        & $4\times 0$ empty matrix &  $\left(\begin{array}{ccccc}
         /
    \end{array}\right)$   \\  \\
    $L_n$  & $\left(\begin{array}{ccccc}
        0  &   0  &   0  &   0  &  0  \\
        0  &   2  &   0  &   0  &  -2  \\
        0  &   0  &   1  &   1  &   0  \\
        0  &   0  &   1  &   2  &   1  \\
        0  &  -2  &   0  &   1  &   3   
    \end{array}\right)$  & $\left(\begin{array}{ccccc}
         2  &   0  &   1  &   -2  \\
         0  &   2  &  -1  &   0  \\
         1  &  -1  &   2  &  -1  \\
        -2  &   0  &  -1  &   2   
    \end{array}\right)$ &    $\left(\begin{array}{ccccc}
         / 
    \end{array}\right)$  \\  \\
    $\beta_{n}$                        & 2               & 1 & 0    \\  \\
    $\text{Spectra}(L_{n})$    & $\{0,0,0.6571,2.5293,4.8136\}$   & $\{0,0.6571,2.5293,4.8136\}$ & / \\ \hline
    \end{tabular}
    \\
    \label{tab:Figure 5 bottom G2}
\end{table}

\begin{table}[H]
    \centering
    \setlength\tabcolsep{6pt}
    \captionsetup{margin=0.9cm}
    \caption{Matrix construction of graph $G_2$ (without isolated points) in the bottom panel of Figure 5.}
    \begin{tabular}{c|ccc}
    \hline
    $n$  & $n=0$ & $n=1$ & $n=2$ \\ \hline 
    $\Omega_n$ & $\text{span}\{e_2,e_3,e_4,e_5\}$ & $\text{span}\{e_{25}, e_{32}, e_{34}, e_{54}\}$ & $\{0\}$ \\ \hline \\
    $B_{n+1}$  & $\begin{array}{@{}r@{}c@{}c@{}c@{}c@{}c@{}l@{}}
            & e_{25} & e_{32} & e_{34} & e_{54}  \\
           \left.\begin{array}{c}
            e_2 \\
            e_3 \\
            e_4 \\
            e_5
            \end{array}\right(
            & \begin{array}{c} -1  \\  0  \\  0  \\  1  \end{array}
            & \begin{array}{c}  1  \\ -1  \\  0  \\  0  \end{array}
            & \begin{array}{c}  0  \\ -1  \\  1  \\  0  \end{array}
            & \begin{array}{c}  0  \\  0  \\  1  \\ -1  \end{array}
            & \left)\begin{array}{c} \\ \\  \\ \\
            \end{array}\right.
        \end{array}$
        & $4\times 0$ empty matrix &  $\left(\begin{array}{ccccc}
         /
    \end{array}\right)$   \\  \\
    $L_n$  & $\left(\begin{array}{ccccc}
        2  &  -1  &   0  &  -1  \\
       -1  &   2  &  -1  &   0  \\
        0  &  -1  &   2  &  -1  \\
       -1  &   0  &  -1  &   2   
    \end{array}\right)$  & $\left(\begin{array}{ccccc}
         2  &  -1  &   0  &  -1  \\
        -1  &   2  &  -1  &   0  \\
         0  &   1  &   2  &   1  \\
        -1  &   0  &   1  &   2   
    \end{array}\right)$ &    $\left(\begin{array}{ccccc}
         / 
    \end{array}\right)$  \\  \\
    $\beta_{n}$                        & 1               & 1 & 0    \\  \\
    $\text{Spectra}(L_{n})$    & $\{0,2,2,4\}$   & $\{0,2,2,4\}$ & / \\ \hline
    \end{tabular}
    \\
    \label{tab:Figure 5 bottom G2}
\end{table}

\begin{table}[H]
    \centering
    \setlength\tabcolsep{6pt}
    \captionsetup{margin=0.9cm}
    \caption{Matrix construction of graph $G_3$ (with isolated points included) in the bottom panel of Figure 5.}
    \begin{tabular}{c|ccc}
    \hline
    $n$  & $n=0$ & $n=1$ & $n=2$ \\ \hline 
    $\Omega_n$ & $\text{span}\{e_1,e_2,e_3,e_4,e_5\}$ & $\text{span}\{e_{25}, e_{32}, e_{34}, e_{54}\}$ & $\{0\}$ \\ \hline \\
    $B_{n+1}$  & $\begin{array}{@{}r@{}c@{}c@{}c@{}c@{}c@{}l@{}}
            & e_{25} & e_{32} & e_{34} & e_{54}  \\
           \left.\begin{array}{c}
            e_1 \\
            e_2 \\
            e_3 \\
            e_4 \\
            e_5
            \end{array}\right(
            & \begin{array}{c}  0 \\ -1  \\  0  \\  0  \\  1  \end{array}
            & \begin{array}{c}  0 \\  1  \\ -1  \\  0  \\  0  \end{array}
            & \begin{array}{c}  0 \\  0  \\ -1  \\  1  \\  0  \end{array}
            & \begin{array}{c}  0 \\  0  \\  0  \\  1  \\ -1  \end{array}
            & \left)\begin{array}{c} \\ \\  \\ \\ \\
            \end{array}\right.
        \end{array}$
        & $4\times 0$ empty matrix &  $\left(\begin{array}{ccccc}
         /
    \end{array}\right)$   \\  \\
    $L_n$  & $\left(\begin{array}{ccccc}
        0  &   0  &   0  &   0  &  0  \\
        0  &   2  &   0  &   0  &  -2  \\
        0  &   0  &   1  &   1  &   0  \\
        0  &   0  &   1  &   2  &   1  \\
        0  &  -2  &   0  &   1  &   3   
    \end{array}\right)$  & $\left(\begin{array}{ccccc}
         2  &   0  &   1  &   -2  \\
         0  &   2  &  -1  &   0  \\
         1  &  -1  &   2  &  -1  \\
        -2  &   0  &  -1  &   2   
    \end{array}\right)$ &    $\left(\begin{array}{ccccc}
         / 
    \end{array}\right)$  \\  \\
    $\beta_{n}$                        & 2               & 1 & 0    \\  \\
    $\text{Spectra}(L_{n})$    & $\{0,0,0.6571,2.5293,4.8136\}$   & $\{0,0.6571,2.5293,4.8136\}$ & / \\ \hline
    \end{tabular}
    \\
    \label{tab:Figure 5 bottom G3}
\end{table}

\begin{table}[H]
    \centering
    \setlength\tabcolsep{6pt}
    \captionsetup{margin=0.9cm}
    \caption{Matrix construction of graph $G_3$ (without isolated points) in the bottom panel of Figure 5.}
    \begin{tabular}{c|ccc}
    \hline
    $n$  & $n=0$ & $n=1$ & $n=2$ \\ \hline 
    $\Omega_n$ & $\text{span}\{e_2,e_3,e_4,e_5\}$ & $\text{span}\{e_{25}, e_{32}, e_{34}, e_{54}\}$ & $\{0\}$ \\ \hline \\
    $B_{n+1}$  & $\begin{array}{@{}r@{}c@{}c@{}c@{}c@{}c@{}l@{}}
            & e_{25} & e_{32} & e_{34} & e_{54}  \\
           \left.\begin{array}{c}
            e_2 \\
            e_3 \\
            e_4 \\
            e_5
            \end{array}\right(
            & \begin{array}{c} -1  \\  0  \\  0  \\  1  \end{array}
            & \begin{array}{c}  1  \\ -1  \\  0  \\  0  \end{array}
            & \begin{array}{c}  0  \\ -1  \\  1  \\  0  \end{array}
            & \begin{array}{c}  0  \\  0  \\  1  \\ -1  \end{array}
            & \left)\begin{array}{c} \\ \\  \\ \\
            \end{array}\right.
        \end{array}$
        & $4\times 0$ empty matrix &  $\left(\begin{array}{ccccc}
         /
    \end{array}\right)$   \\  \\
    $L_n$  & $\left(\begin{array}{ccccc}
        2  &  -1  &   0  &  -1  \\
       -1  &   2  &  -1  &   0  \\
        0  &  -1  &   2  &  -1  \\
       -1  &   0  &  -1  &   2   
    \end{array}\right)$  & $\left(\begin{array}{ccccc}
         2  &  -1  &   0  &  -1  \\
        -1  &   2  &  -1  &   0  \\
         0  &   1  &   2  &   1  \\
        -1  &   0  &   1  &   2   
    \end{array}\right)$ &    $\left(\begin{array}{ccccc}
         / 
    \end{array}\right)$  \\  \\
    $\beta_{n}$                        & 1               & 1 & 0    \\  \\
    $\text{Spectra}(L_{n})$    & $\{0,2,2,4\}$   & $\{0,2,2,4\}$ & / \\ \hline
    \end{tabular}
    \\
    \label{tab:Figure 5 bottom G3}
\end{table}

\begin{table}[H]
    \centering
    \setlength\tabcolsep{6pt}
    \captionsetup{margin=0.9cm}
    \caption{Matrix construction of graph $G_2$ in the bottom panel of Figure 5.}
    \begin{tabular}{c|ccc}
    \hline
    $n$  & $n=0$ & $n=1$ & $n=2$ \\ \hline 
    $\Omega_n$ & $\text{span}\{e_1,e_2,e_3,e_4\}$ & $\text{span}\{e_{12}, e_{14}, e_{32}, e_{34}\}$ & $\{0\}$ \\ \hline \\
    $B_{n+1}$  & $\begin{array}{@{}r@{}c@{}c@{}c@{}c@{}c@{}l@{}}
            & e_{12} & e_{14} & e_{32} & e_{34} & e_{54}  \\
           \left.\begin{array}{c}
            e_1 \\
            e_2 \\
            e_3 \\
            e_4 \\
            e_5
            \end{array}\right(
            & \begin{array}{c}  0 \\  0  \\  0  \\  0  \\  0  \end{array}
            & \begin{array}{c}  0 \\ -1  \\  0  \\  0  \\  1  \end{array}
            & \begin{array}{c}  0 \\  0  \\ -1  \\ -1  \\  0  \end{array}
            & \begin{array}{c}  0 \\  0  \\  0  \\  1  \\  1  \end{array}
            & \begin{array}{c}  0 \\  1  \\  0  \\  0  \\ -1  \end{array}
            & \left)\begin{array}{c} \\ \\  \\ \\ \\
            \end{array}\right.
        \end{array}$
        & $4\times 0$ empty matrix &  $\left(\begin{array}{ccccc}
         /
    \end{array}\right)$   \\  \\
    $L_n$  & $\left(\begin{array}{ccccc}
         2  &  -1  &   0  &  -1  \\
        -1  &   2  &  -1  &   0  \\
         0  &  -1  &   2  &  -1  \\
        -1  &   0  &  -1  &   2   
    \end{array}\right)$  & $\left(\begin{array}{ccccc}
         2  &   1  &   1  &   0  \\
         1  &   2  &   0  &   1  \\
         1  &   0  &   2  &   1  \\
         0  &   1  &   1  &   2   
    \end{array}\right)$ &    $\left(\begin{array}{ccccc}
         / 
    \end{array}\right)$  \\  \\
    $\beta_{n}$                        & 1               & 1 & 0    \\  \\
    $\text{Spectra}(L_{n})$    & $\{0,2,2,4\}$   & $\{0,2,4,4\}$ & / \\ \hline
    \end{tabular}
    \\
    \label{tab:Figure 5 bottom G2}
\end{table}

\begin{table}[H]
    \centering
    \setlength\tabcolsep{6pt}
    \captionsetup{margin=0.9cm}
    \caption{Matrix construction of graph $G_4$ in the bottom panel of Figure 5.}
    \begin{tabular}{c|ccc}
    \hline
    $n$  & $n=0$ & $n=1$ & $n=2$ \\ \hline 
    $\Omega_n$ & $\text{span}\{e_1,e_2,e_3,e_4,e_5\}$ & $\text{span}\{e_{13}, e_{25}, e_{32}, e_{34},e_{45}\}$ & $\{0\}$ \\ \hline \\
    $B_{n+1}$  & $\begin{array}{@{}r@{}c@{}c@{}c@{}c@{}c@{}l@{}}
            & e_{13} & e_{25} & e_{32} & e_{34} & e_{45}  \\
           \left.\begin{array}{c}
            e_1 \\
            e_2 \\
            e_3 \\
            e_4 \\
            e_5
            \end{array}\right(
            & \begin{array}{c}-1 \\  0 \\  1  \\  0  \\  0  \end{array}
            & \begin{array}{c} 0 \\ -1 \\  0  \\  0  \\  1  \end{array}
            & \begin{array}{c} 0 \\  1 \\ -1  \\  0  \\  0  \end{array}
            & \begin{array}{c} 0 \\  0 \\ -1  \\  1  \\  0  \end{array}
            & \begin{array}{c} 0 \\  0 \\  0  \\  1  \\ -1  \end{array}
            & \left)\begin{array}{c} \\ \\  \\ \\ \\
            \end{array}\right.
        \end{array}$
        & $5\times 0$ empty matrix &  $\left(\begin{array}{ccccc}
         /
    \end{array}\right)$   \\  \\
    $L_n$  & $\left(\begin{array}{ccccc}
         1  &   0  &  -1  &   0  &   0  \\
         0  &   2  &  -1  &   0  &  -1  \\
        -1  &  -1  &   3  &  -1  &   0  \\
         0  &   0  &  -1  &   2  &  -1  \\
         0  &  -1  &   0  &  -1  &   2   
    \end{array}\right)$  & $\left(\begin{array}{ccccc}
         2  &   0  &  -1  &  -1  &   0  \\
         0  &   2  &  -1  &   0  &  -1  \\
        -1  &  -1  &   2  &   1  &   0  \\
        -1  &   0  &   1  &   2  &   1  \\
         0  &  -1  &   0  &   1  &   2   
    \end{array}\right)$ &    $\left(\begin{array}{ccccc}
         / 
    \end{array}\right)$  \\  \\
    $\beta_{n}$                        & 1               & 1 & 0    \\  \\
    $\text{Spectra}(L_{n})$    & $\{0,0.8299,2,2.6889,4.4812\}$   & $\{0,0.8299,2,2.6889,4.4812\}$ & / \\ \hline
    \end{tabular}
    \\
    \label{tab:Figure 5 bottom G4}
\end{table}

\begin{landscape}
\begin{table}[H]\footnotesize
    \centering
    \setlength\tabcolsep{6pt}
    \captionsetup{margin=0.9cm}
    \caption{Matrix construction of graph $G_5$ in the bottom panel of Figure 5.}
    \begin{tabular}{c|ccc}
    \hline
    $n$  & $n=0$ & $n=1$ & $n=2$ \\ \hline 
    $\Omega_n$ & $\text{span}\{e_1,e_2,e_3,e_4,e_5\}$ & $\text{span}\{e_{12}, e_{13}, e_{14}, e_{15}, e_{25}, e_{32}, e_{34}, e_{54}\}$ & $\text{span}\{e_{125},e_{132},e_{134},e_{154}\}$ \\ \hline \\
    $B_{n+1}$  & $\begin{array}{@{}r@{}c@{}c@{}c@{}c@{}c@{}c@{}c@{}c@{}l@{}}
            & e_{12} & e_{13} & e_{14} & e_{15} & e_{25} & e_{32} & e_{34} & e_{54}  \\
           \left.\begin{array}{c}
            e_1 \\
            e_2 \\
            e_3 \\
            e_4 \\
            e_5
            \end{array}\right(
            & \begin{array}{c} -1 \\  1  \\  0  \\  0  \\ 0  \end{array}
            & \begin{array}{c} -1 \\  0  \\  1  \\  0  \\ 0  \end{array}
            & \begin{array}{c} -1 \\  0  \\  0  \\  1  \\ 0  \end{array}
            & \begin{array}{c} -1 \\  0  \\  0  \\  0  \\ 1  \end{array}
            & \begin{array}{c}  0 \\ -1  \\  0  \\  0  \\ 1  \end{array}
            & \begin{array}{c}  0 \\  1  \\ -1  \\  0  \\ 0  \end{array}
            & \begin{array}{c}  0 \\  0  \\ -1  \\  1  \\ 0  \end{array}
            & \begin{array}{c}  0 \\  0  \\  0  \\  1  \\-1  \end{array}
            & \left)\begin{array}{c} \\ \\  \\ \\ \\ 
            \end{array}\right.
        \end{array}$
        & $\begin{array}{@{}r@{}c@{}c@{}c@{}c@{}c@{}c@{}c@{}l@{}}
            & e_{125} & e_{132} & e_{134} & e_{154} \\
           \left.\begin{array}{c}
            e_{12} \\
            e_{13} \\
            e_{14} \\
            e_{15} \\
            e_{25} \\
            e_{32} \\
            e_{34} \\
            e_{54}
            \end{array}\right(
            & \begin{array}{c}  1 \\  0  \\  0  \\ -1  \\ 1 \\ 0 \\ 0  \\ 0   \end{array}
            & \begin{array}{c} -1 \\  1  \\  0  \\  0  \\ 0 \\ 1 \\ 0  \\ 0   \end{array}
            & \begin{array}{c}  0 \\  1  \\ -1  \\  0  \\ 0 \\ 0 \\ 1  \\ 0   \end{array}
            & \begin{array}{c}  0 \\  0  \\ -1  \\  1  \\ 0 \\ 0 \\ 0  \\ 1   \end{array}
            & \left)\begin{array}{c} \\ \\  \\ \\ \\ \\ \\ \\
            \end{array}\right.
        \end{array}$ &  $4\times 0$ empty matrix  \\  \\
    $L_n$  & $\left(\begin{array}{ccccc}
         4  &  -1  &  -1  &  -1  &  -1  \\
        -1  &   3  &  -1  &   0  & -1  \\
        -1  &  -1  &   3  &  -1  &  0  \\
        -1  &   0  &  -1  &   3  & -1  \\
        -1  &  -1  &   0  &  -1  &  3  
    \end{array}\right)$  & $\left(\begin{array}{cccccccccc}
         4  &   0  &   1  &   0  &  0  &  0  &  0  &  0 \\
         0  &   4  &   0  &   1  &  0  &  0  &  0  &  0 \\
         1  &   0  &   4  &   0  &  0  &  0  &  0  &  0 \\
         0  &   1  &   0  &   4  &  0  &  0  &  0  &  0 \\
         0  &   0  &   0  &   0  &  3  & -1  &  0  & -1 \\
         0  &   0  &   0  &   0  & -1  &  3  &  1  &  0 \\
         0  &   0  &   0  &   0  &  0  &  1  &  3  &  1 \\
         0  &   0  &   0  &   0  & -1  &  0  &  1  &  3 \\
    \end{array}\right)$ &    $\left(\begin{array}{ccccc}
         3 & -1 & 0 & -1 \\
        -1 &  3 & 1 &  0 \\
         0 &  1 & 3 &  1 \\
        -1 &  0 & 1 &  3
    \end{array}\right)$  \\  \\
    $\beta_{n}$                        & 1               & 0 & 0    \\  \\
    $\text{Spectra}(L_{n})$    & $\{0,3,3,5,5\}$   & $\{1,3,3,3,3,5,5,5\}$ & $\{1,3,3,5\}$ \\ \hline
    \end{tabular}
    \\
    \label{tab:Figure 5 bottom G5}
\end{table}

\end{landscape}

\bibliographystyle{unsrt}
\bibliography{refs}
\end{document}